\documentclass[review,11pt]{elsarticle}

\usepackage{lineno,hyperref}

\journal{European Journal of Operational Research}

\newcommand{\halmos}{\ensuremath{\Box}}
\newcommand{\p}{\phantom}
\newcommand{\braces}[1]{\left\{ #1 \right \}}

\newcommand{\lb}{\linebreak}
\newcommand{\comment}[1]{}
\renewcommand{\mid}{\,|\,}

\usepackage{amssymb}
\usepackage{amsthm}
\usepackage{lipsum} 
\usepackage{amsmath}
\usepackage{tabularx}
\usepackage{bbm}
\usepackage{booktabs}
\usepackage{tikz}
\usetikzlibrary{shapes}
\modulolinenumbers[5]
\usepackage{enumitem}
\setlist{align=left,topsep=0pt,itemsep=-1ex,partopsep=1ex,parsep=1ex}
\setlist[enumerate]{topsep=0pt,leftmargin=*,labelsep=2ex,itemsep=-1ex}

\usepackage[margin=1.5cm]{geometry}

\selectfont

\newtheorem{theorem}{Theorem}
\newtheorem{corollary}{Corollary}
\newtheorem{proposition}{Proposition}

\newdefinition{axiom}{Axiom}
\newdefinition{definition}{Definition}

\definecolor{blue}{rgb}{0.23,0.58,0.89}

\biboptions{authoryear}

\begin{document}

\begin{frontmatter}

\title{Decision Programming for Mixed-Integer Multi-Stage Optimization under Uncertainty}


\author{Ahti Salo\corref{cor}}
\ead{ahti.salo@aalto.fi}

\author{Juho Andelmin}
\author{Fabricio Oliveira}

\cortext[cor]{Corresponding author}
\address{Systems Analysis Laboratory, Department of Mathematics and Systems Analysis \\ Aalto University School of Science, P.O. Box 11100, 00076 Aalto, FINLAND}

\begin{abstract}
Influence diagrams are widely employed to represent multi-stage decision problems in which each decision is a choice from a discrete set of alternatives, uncertain chance events have discrete outcomes, and prior decisions may influence the probability distributions of uncertain chance events endogenously. In this paper, we develop the \textit{Decision Programming} framework which extends the applicability of influence diagrams by developing mixed-integer linear programming formulations for solving such problems in the presence of many kinds of constraints. In particular, Decision Programming makes it possible to (i) solve problems in which earlier decisions cannot necessarily be recalled later, for instance, when decisions are taken by agents who cannot communicate with each other; (ii) accommodate a broad range of deterministic and chance constraints, including those based on resource consumption, logical dependencies or risk measures such as Conditional Value-at-Risk; and (iii) determine all non-dominated decision strategies in problems which involve multiple value objectives. In project portfolio selection problems, Decision Programming allows scenario probabilities to depend endogenously on project decisions and can thus be viewed as a generalization of \textit{Contingent Portfolio Programming}.  %(\citealp{GustafssonSalo2005}) 
We provide several illustrative examples as well as evidence on the computational performance of Decision Programming formulations.
\end{abstract}


\end{frontmatter}

\clearpage

%\linenumbers

\section{Introduction}

Influence diagrams, in their many variants (see, e.g., \citealp{Bielza2011, Diehl2004, DiezEtAl2018,HowardMatheson1984,HowardMatheson2005}), are widely employed to represent decision problems whose consequences depend on discrete sets of both uncertain chance events and decisions which are taken in multiple stages. Specifically, such decisions and chance events are represented by decision and chance nodes in an acyclic graph whose arcs indicate (i) what information is available to the decision maker (DM) and (ii) how realizations of chance events depend on earlier decisions and chance events. The value node represents consequences which are associated with the DM's decisions and the realization of chance events. 

The DM's risk preferences are typically modeled with a utility function over the set of consequences. The optimal solution to the influence diagram is the decision strategy that, at each decision node, assigns a decision alternative to every possible state of information at this node so that the combination of these decisions maximizes the DM's expected utility. If the diagram fulfills the \emph{`no forgetting' assumption}, meaning that earlier decisions can be recalled when making later ones (see, e.g., \citealp{LauritzenNilsson2001,JorgensenEtAl2014}), this optimal strategy can be computed with well-established techniques, for example by carrying out local transformations such as arc reversals and node removals \citep{Shachter1986, Shachter1988}, or by formulating the equivalent decision tree representation and solving it with dynamic programming \citep{TatmanShachter1990}. While this assumption often holds, there are problems in which it does not, for example, in distributed decision making by agents such as military patrols who may not be able to communicate with each other (for examples, see, e.g., \citealp{Zhang1994}). Examples include problems of adversarial risk analysis in which agents may not be able to observe what decisions the other agents have taken \citep{RiosInsua2009, Roponen2020}. Another important context is the risk management of safety-critical systems that must not fail even if information cannot be synchronized due to disruptions or communication delays.  
More generally, dynamic programming is a restrictive solution approach, because the optimal strategy within a branch that unfolds from a given decision node cannot depend on decisions in other branches of the decision tree. Thus, the objective function cannot include risk measures such as Value-at-Risk, which reflects the full variability of consequences across the entire decision tree. Project portfolio selection problems, too, give rise to analogous dependencies, because the consumption of shared resources, for example, implies that optimal decisions for one project are contingent on those for others. Consequently, the optimal strategy for a given project cannot be determined without considering strategies for the other projects \citep{GustafssonSalo2005}. 

In this paper, we develop the \emph{Decision Programming} modeling framework which (i) uses influence diagrams to represent the structure of discrete multi-stage decision problems under uncertainty, including those that cannot be solved with dynamic programming techniques, (ii) extends the problem formulation by allowing for the specification of both  deterministic (e.g., logical dependencies, costs arising at one or more nodes) and chance (i.e., probabilistic) constraints, and (iii) converts the resulting problem representations into mixed-integer linear programming (MILP) problems. In particular, this framework is generic enough for solving \underline{li}mited \underline{m}emory \underline{i}nfluence \underline{d}iagrams (LIMID) in which the `no forgetting' assumption may not hold. Decision programming also makes it possible to compute all non-dominated strategies in problems with  multiple objectives, represented by multiple value nodes. Importantly, all these modeling features can be incorporated into corresponding MILP problems that can be solved efficiently with available software tools (for a survey, see, e.g., \citealp{Fourer2017}), as evidenced by our computational results which demonstrate that problems of considerable size can be solved to optimality.  

Our contribution is relevant to Stochastic Programming \citep{BirgeLouveaux2011} as it provides a general framework for problems in which decisions are made over several stages and realizations of uncertain events are observed between pairs of successive stages. In the first stage, an initial decision is selected, and subsequent recourse decisions are selected after observing the realizations of uncertain earlier events. We distinguish between endogenous and exogenous uncertainties based on whether earlier decisions can influence conditional probability distributions. Both types of uncertainties can be accommodated in Decision Programming by converting influence diagrams and adjoining constraints into multi-stage stochastic integer programming (MSSIP) problems that can be efficiently solved using off-the-shelf MILP solvers. That is, the diagram is first converted into a sequence of decision and chance nodes. This sequence is then employed when transforming the deterministic equivalent MILP formulation of the MSSIP. 

The rest of this paper is structured as follows. Section 2 discusses earlier approaches. Section 3 develops the Decision Programming framework. Section 4 presents illustrative examples. Section 5 develops modeling approaches for dealing with risk preferences, chance constraints and multiple objectives. Section 6 gives results on computational performance. Section 7 concludes and provides directions for further development of the framework. 

\section{Earlier approaches}

Influence diagrams were initially developed in the 1970's (\citealp{Olmsted1983, HowardMatheson1984,HowardMatheson2005,HowardMatheson2006,HowardEtAl2006}) to represent informational and probabilistic dependencies between decisions and uncertain chance events which, taken together, determine consequences for the DM. If the `no-forgetting' assumption holds so that earlier decisions are known when making later ones, and the aim is to maximize expected utility at the value node, these diagrams can be solved with well-established techniques, for instance, by forming an equivalent decision tree that can be solved through dynamic programming \citep{TatmanShachter1990}; or by removing decision and chance nodes from the diagram one-by-one, possibly after arc reversals (see, e.g., \citealp{Shachter1986,SmithEtAl1993,HowardMatheson2005}). 

As visual tools for problem representation, influence diagrams differ from decision trees in suggesting that the problem has a symmetric structure in which the sets of  alternative decisions as well as realizations of chance events do not depend on preceding decision and chance nodes. Still, asymmetric problems can be modeled with influence diagrams through the appropriate definition of node states and their dependencies (see, e.g., \citep{SmithEtAl1993}). Mathematically, the mapping of input parameters (i.e., probabilities, decisions) to outputs (i.e., expected utilities) in influence diagrams is a piecewise multilinear function \citep{BorgonovoTonoli2014}, which fact underpins developments in this paper. For an account of the evolution of influence diagrams, see \citet{Bielza2011}. 

Problems in which earlier decisions cannot be recalled give rise to LIMIDs which are computationally challenging because optimal strategies cannot be determined through an equally straightforward series of local computations. \citet{Zhang1994} discuss these and other kinds of influence diagrams. \citet{LauritzenNilsson2001} develop an iterative policy updating approach for LIMIDs by solving a series of expected utility maximization problems by message passing in a junction tree derived from the influence diagram. \citet{HovgaardBrincker2016} describe an application of LIMIDs to structural damage protection. \citet{Maua2016} study the computational performance of $k$-neighborhood local search algorithms and propose approximate algorithms. Further optimization formulations for solving influence diagrams in junction trees are presented by \citet{ParmentierEtAl2019}.

However, these approaches based on local computations and iterative message passing schemes are ill-equipped for problems involving constraints that span across the entire problem  (e.g., due to logical interdependencies, limited budgets, bounds on risk levels) and whose fulfilment cannot be determined locally. For example, the DM may seek to maximize the expected net present value (NPV) subject to the requirement that the expectation in the lower tail of the NPV distribution is not too low (i.e., Conditional Value-at-Risk, which is a coherent risk measure; \citealp{ArtznerEtAl1999}). These local approaches also encounter difficulties in problems in which several objectives associated with their respective multiple value nodes have to be explicitly addressed (see, e.g., \citealp{Diehl2004}). 

In portfolio decision analysis \citep{SaloEtAl2011}, influence diagrams help portray the overall structure of probabilistic and informational dependencies, but they cannot handle constraints arising from limited budgets or logical dependencies between alternatives. For project selection problems, \textit{Contingent Portfolio Programming} \citep{GustafssonSalo2005} employs MILP to determine optimal project management strategies when the projects' cash flows are contingent on scenarios whose probabilities cannot depend endogenously on project decisions. \citet{Vilkkumaa2018} extend this approach to single-stage selection problems in which scenario probabilities can depend endogenously on project decisions. \citet{LiesioSalo2012} derive decision recommendations for single-stage project selection problems with one objective and possibly incomplete utility and probability information. Yet, none of these earlier approaches are equipped to handle problems in which there is a {\it combination} of endogenous uncertainties, several decision stages, and multiple objectives.

Several papers use stochastic programming as the underpinning framework for modeling multi-stage problems under uncertainty. Nevertheless, the literature on endogenous uncertainty in stochastic programming is still sparse, because the existing models depart from domains in which well performing solution techniques are available, most prominently convex programming in general, and linear programming in particular. 

Most of the stochastic programming literature focuses on problems in which decisions can influence the information structure, in particular the timing of unveiling  uncertainties, as opposed to the actual probability distributions associated with uncertain events. \citet{goel2006} develop a stochastic programming formulation for multi-stage problems for the timing of oil well exploitation, which is assumed not to influence the uncertain amount of recoverable oil. Building on developments in \citet{goel2004}, they propose a unified framework and solution methods to handle problems in which the decisions influence the time of observing uncertainties. \citet{gupta2011, gupta2014} present specialized solution methods for oil and gas field development. \citet{colvin2008} propose a stochastic programming model for novel product development in pharmaceutical research, further extended by \citet{colvin2009}. In this context, the timing of when uncertainties are resolved is influenced endogenously by the decisions on how to perform clinical trials which, however, leads to computational challenges \citep{colvin2010}. \citet{solak2010} address R\&D project portfolio optimization under endogenous uncertainty, acknowledging that the inclusion of decision dependent uncertainties significantly degrades tractability. To tackle this issue, they propose a sophisticated solution method, exploiting the formulation devised specifically for the problem. \citet{apap2017} provide a comprehensive recent literature overview and propose an approach for problems with a decision-dependent information structure. 

Problems where decisions can (also) affect the probability distributions of uncertain events have been much less explored. The predominant strategy has been to remove decision dependent probabilities using appropriate transformations in the probability measure, as described by \citet{rubinstein1993} (see also \citealp{pflug2012}), or in the probability distribution itself (cf. \citealp{dupacova2006}). In their overview of this scarce literature, \citet{hellemo2018} propose a taxonomy of distinct classes for stochastic programs with endogenous uncertainties and possible formulation approaches. They also report computational experiments to highlight how challenging these problems are for state-of-the-art optimization solvers. 

In fact, multi-stage optimization problems under uncertainty can involve decision dependent probabilities, parameters, and/or information structures \citep{hellemo2018}. The \emph{Decision Programming} framework developed here is general enough to encompass all these variants, on condition that each chance event has a finite number of possible realizations and decisions correspond to choices from a finite set of discrete alternatives; this is the case in the majority of the aforementioned approaches and the following development. 

\section{Methodological development} 

\subsection{Influence diagram representation of the decision problem}

Multi-stage decision problems under uncertainty can be modeled as acyclic networks $G = (N,A)$ whose nodes $N = C \cup D \cup V$ consist of chance nodes $C$, decision nodes $D$, and value nodes $V$. Chance nodes $C$ represent uncertain events associated with random variables; decision nodes $D$ correspond to decisions among discrete alternatives; and value nodes $V$ represent consequences that result from the realizations of random variables at chance nodes and the decisions made at decision nodes.  

Dependencies between nodes are represented by arcs $A = \{ (i,j) ~|~ i,j \in N \}$. A \textit{path} of length $k$ is a sequence of nodes $(i_1, i_2, \ldots, i_k)$ such that $(i_l,i_{l+1}) \in A$ for all $l = 1,\ldots, k-1$. The \emph{information set} of a node $j \in N$, defined as $I(j) = \{i \in N ~|~ (i,j) \in A \}$, consists of the direct predecessors of $j$ from which there is an arc to $j$. Since the network $G$ is acyclic, the nodes $N$ can be indexed consecutively with integers $1, 2, \dots, |N|$ (where $|\,\cdot\,|$ denotes the number of elements in a set) so that for each node $j \in N$, the indices of the nodes in its information set $I(j)$ are smaller than $j$ (i.e., $i < j$ for all $i \in I(j)$).

We denote the number of chance nodes by $n_C = |C|$ and the number of decision nodes by $n_D = |D|$. These $n = n_C + n_D$ chance and decision nodes are indexed as $C \cup D = \{ 1,2,\dots, n \}$, while the $n_V = |V| = |N| - n $ value nodes are indexed as $n + 1,\ldots, n + n_V$. For now, we assume that there is a single value node in the influence diagram (the extension to multiple value nodes is covered in Section \ref{sec:multiple}). Consequences at this value node are determined by the decisions and the realization of chance events which do not depend on consequences. Thus, there are no arcs $(i,j) \in A$ such that $i \in V$ and $j \in C \cup D$.

Each chance and decision node $j \in C \cup D$ has a finite set $S_j$ of discrete states. The occurrence of states depend on their possible {\it information states} $s_{I(j)} \in S_{I(j)} = \prod_{i \in I(j)} S_i$, defined as combinations of states of all nodes in the information set $I(j)$. For each chance node $j \in C$, these states correspond to realizations of the random variable $X_j$, which depends probabilistically on the states $s_i$ of the nodes $i \in I(j)$ in the information set of $j$. For a decision node $j \in D$, each state $s_j \in S_j$ corresponds to a decision that is made based on the information state $s_{I(j)}$. For brevity, we use $X_j, j = 1,\ldots,n$, to denote both random variables which are associated with chance nodes $j \in C$ and decision variables which are associated with decision nodes $j \in D$. 

Specifically, if $j \in C$ is a chance node whose information state is $s_{I(j)}$, then state $s_j \in S_j$ occurs with the conditional probability 
\begin{equation}\label{CondProb}
\mathbb{P}(X_j = s_j \mid X_{I(j)} = s_{I(j)}), \qquad \forall \,j \in C, \, s_j \in S_j, \, s_{I(j)} \in S_{I(j)},  
\end{equation} 
where $X_{I(j)} = s_{I(j)}$ means that the states of the variables $X_i$ in the information set $i \in I(j)$ are the same as specified by the information state $s_{I(j)}$. For each decision node $j \in D$, a {\it local decision strategy} $Z_j : S_{I(j)} \mapsto S_j$ is a function that maps each information state in $S_{I(j)}$ to a decision in $S_j$. A \textit{(global) decision strategy} $Z$ is a set of local decision strategies which contains one local strategy $Z_j$ for each decision node $j \in D$. The set of all decision strategies is denoted by $\mathbb{Z}$.  

\subsection{Paths}

A \textit{path} $s = (s_1,s_2,\ldots,s_n)$ of length $n$ is a sequence of states $s_i \in S_i$ of all chance and decision nodes, i.e., $i \in C\cup D$ for all $i=1,\ldots,n$. The set $S$ of all paths of length $n$ is 
\begin{equation}\label{eq:SPATH}
S = S_{1:n} = \{ (s_1,s_2,\ldots,s_n) \mid s_i \in S_i,\ i=1,\ldots,n \}.
\end{equation} 
Paths of length $k < n$ are sequences $s_{1:k} = (s_1,s_2,\ldots,s_k)$ such that $s_i \in S_i, ~i \leq k$. If $s_{1:k} \in S_{1:k}, k < n$, and $s_{k+1} \in S_{k+1}$, the state $s_{k+1}$ can be appended to $s_{1:k}$ to form the path $s_{1:k+1} = (s_1,s_2,\dots,s_k,s_{k+1}) \in S_{1:k+1}$. If $s_{1:k} \in S_{1:k}, ~k \leq n$, and $I \subsetneq \{1 , \ldots, k \}$, then $s_{I}$ is a  subsequence of $s_{1:k}$ for the nodes $i \in I$. Thus, $s_I$ is a sequence of length $|I|$ which contains the same states as $s_{1:k}$ for nodes $i \in I$. 

A decision strategy $Z \in \mathbb{Z}$ is \emph{compatible} with the path $s \in S$ if and only if $ Z_j(s_{I(j)}) = s_j, \ \forall \, Z_j \in Z, j \in D$. Thus, at each decision node $j \in D$, $Z_j \in Z$ maps the information state $s_{I(j)}$ contained in $s$ to the corresponding decision $s_j$ in $s$. In this case, the probability of path $s$ is $\mathbb{P}(s \mid Z) = \prod_{i \in C} \mathbb{P}(X_i= s_i \mid X_{I(i)} = s_{I(i)})$. On the other hand, if $Z$ is not compatible with $s$, it contains some local decision strategy $Z_j$, $j \in D$, such that the information state $s_{I(j)}$ contained in $s$ is mapped to a decision which differs from the state $s_j$ of node $j$ in $s$. As a result, choosing $Z$ means that $s$ cannot occur and therefore $\mathbb{P}(s \mid Z) = 0$. 

More generally, for a given decision strategy $Z \in \mathbb{Z}$, the probability of a path $s \in S$ can be expressed recursively as a function of the conditional probabilities (\ref{CondProb}) and local decision strategies so that  
\begin{equation} \label{path_probabilities}
\mathbb{P}(s_{1:k} \mid Z)= \bigg( \prod_{\stackrel{i \in C}{i \leq k}} \mathbb{P} \big(X_i=s_i \mid X_{I(i)} = s_{I(i)} \big) \bigg) \bigg( \prod_{\stackrel{j \in D}{j \leq k}} \mathbb{I} \big( Z_j(s_{I(j)})=s_j \big) \bigg), 
\end{equation}
where the indicator function $\mathbb{I}(\,\cdot\,)$ is defined so that 
\begin{equation}
    \mathbb{I}(Z_j(s_{I(j)}) = s_j) = \begin{cases}1, &\text{if } Z_j(s_{I(j)})=s_j, \\
                                       0, &\text{otherwise.} \end{cases}
\end{equation}

\subsection{Characterizing path probabilities using linear inequalities}

A given decision strategy $Z\in \mathbb{Z}$ assigns probabilities to all paths $s_{1:k} \in S_{1:k}$, $k=1,\ldots,n$, in accordance with \eqref{path_probabilities}. However, this expression does not suggest efficient ways of computing these probabilities. One could introduce binary variables taking values of the indicator functions $\mathbb{I} \big(Z_j(s_{I(j)})=s_j \big),$ for all $j \in D,$ whose multiplication would lead to a mixed-integer nonlinear programming (MINLP) problem which could be converted into a equivalent MILP. An early version of the Decision Programming approach relied on this strategy, which, despite being feasible, led to a MILP formulation with a weak linear programming relaxation, and was therefore deemed too inefficient for off-the-shelf solver performance.  

Alternatively, we characterize the probabilities of paths $s_{1:k} \in S_{1:k}$, $k = 1,\dots,n$, through sets of linear inequalities. Towards this end, local decision strategies $Z_j$, $j \in D$, are modelled through corresponding binary variables $z(s_j \mid s_{I(j)}) \in \{0,1\}$ such that $z(s_j \mid s_{I(j)}) = 1$ if and only if $Z_j$ maps the information state $s_{I(j)}$ to the decision $s_j \in S_j$, i.e.,
\begin{equation}\label{zvariables}
Z_j(s_{I(j)}) = s_j \iff z(s_j \mid s_{I(j)}) = 1, \qquad \forall \,j \in D,\, s_j  \in S_j, \,s_{I(j)} \in S_{I(j)}. 
\end{equation} 
Furthermore, the mutual exclusivity of the decisions is ensured through the constraints
\begin{equation}\label{normconstraint}
\sum_{s_j \in S_j} z(s_j \mid s_{I(j)}) = 1, \qquad  \forall \,j \in D,  \,s_{I(j)} \in S_{I(j)}, 
\end{equation}
which ensure that exactly one decision $s_j \in S_j$ is chosen for every information state $s_{I(j)} \in S_{I(j)}$. 

For the given decision strategy $Z \in \mathbb{Z}$, the corresponding probability $\pi(s)$ of any path $s \in S$ can be derived recursively as follows. To initialize the recursive process, let $\pi_0(s) = 1$. Suppose that the probabilities $\pi_i(s) = \mathbb{P}(X_{1:k-1} = s_{1:k-1} \mid Z) $ are known for nodes $i \leq k-1$ and consider the next node $k \leq n$. If $k \in C$ is a chance node, let
\begin{equation}\label{condprob}
\pi_k(s) = \mathbb{P}\left(X_k = s_k \mid X_{I(k)} = s_{I(k)}\right) \pi_{k-1}(s), 
\end{equation}
where the first term on the right side of \eqref{condprob} is given by \eqref{CondProb}. If $k \in D$ is a decision node, let  
\begin{equation} \label{DecisionPi}
\pi_k(s) = \begin{cases}
\pi_{k-1}(s), & \quad \text{if } z(s_k \mid s_{I(k)}) = 1 \\ 
0, & \quad \text{if } z(s_k \mid s_{I(k)}) = 0. \end{cases}  
\end{equation}
This assignment corresponds to the inequalities 
\begin{equation*}\label{decisionconstraints}
\text{max} \{0, \, \pi_{k-1}(s) + z(s_k \mid s_{I(k)}\textbf{}) - 1\}  \leq \pi_k(s) \leq \text{min} \{ \pi_{k-1}(s),\, z(s_k \mid s_{I(k)}) \}, 
\end{equation*}
which are equivalent to 
\begin{eqnarray}
\pi_k(s) & \leq & \pi_{k-1}(s) \label{inequality1} \\
\pi_k(s) & \leq & z(s_k \mid s_{I(k)}) \label{inequality2} \\
\pi_k(s) & \geq & 0 \label{inequality3} \\
\pi_k(s) & \geq & \pi_{k-1}(s) + z(s_k \mid s_{I(k)}) - 1. \label{inequality4}
\end{eqnarray}

Theorem \ref{scenariopaths} states that the path probabilities implied by any strategy $Z$ can be calculated through the assignment \eqref{zvariables}--\eqref{DecisionPi}. Importantly, the equivalence between the assignments \eqref{zvariables}--\eqref{DecisionPi} and the inequalities \eqref{inequality1}--\eqref{inequality4} implies that the strategy which maximizes the expectation of a real-valued function over paths can be determined through optimization by employing these inequalities as constraints on the decision variables $z(s_k|s_{I(k)}), k \in D, s_k \in S_k, s_{I(k)} \in S_{I(k)}$.

\begin{theorem} \label{scenariopaths}
Let $Z \in \mathbb{Z}$ be a decision strategy and choose a path $s \in S$. If $\pi_k(s)$, $k = 1, \dots, n$, and $z(s_j \mid s_{I(j)})$, $\forall j \in D$, satisfy the constraints \eqref{zvariables}--\eqref{DecisionPi}, then 
\begin{equation} \label{pi_recursion}
\pi_k(s) = \mathbb{P}(X_{1:k} = s_{1:k} \mid Z), \qquad \forall \, k = 1,\ldots,n. \end{equation}
In particular, $\pi(s) \overset{def}{=} \pi_n(s)$ is the probability of the path $s$ for the decision strategy $Z$.
\end{theorem}

{\noindent \bf Proof}. See Appendix A. 

\subsection{Maximization of expected utility}

We assume that at the value node $v \in V$, the function $Y_v : S_{I(v)} \mapsto \mathbb{C}$ maps combinations of states of the nodes in its information set $I(v)$ to a set of consequences $\mathbb{C}$ and that there exists a real-valued utility function $U : \mathbb{C} \mapsto \mathbb{R}$ that is defined over $\mathbb{C}$. Then, the utility associated with the path $s \in S$ can be precomputed as
\begin{equation}\label{preutil}
\mathcal{U}(s) = U[Y_v(s_{I(v)})].   
\end{equation}

Because the path probabilities $\pi(s)$, $s \in S$, for the selected decision strategy $Z \in \mathbb{Z}$ are given by Theorem \ref{scenariopaths}, it follows that the decision strategy which maximizes the DM's expected utility is the solution to the optimization problem in Corollary \ref{cor:opt1}. 

\begin{corollary}\label{cor:opt1}
The expected utility is maximized by the decision strategy $Z\in \mathbb{Z}$ which solves the optimization problem  
\begin{equation}\label{optimization}
\mathop{\emph{maximize}}\limits_{Z \in \mathbb{Z}} \sum_{s \in S} \pi(s) \mathcal{U}(s)
\end{equation}
subject to constraints \eqref{zvariables}--\eqref{condprob} and \eqref{inequality1}--\eqref{inequality4} on decision variables $z(s_k \mid s_{I(k)}) \in \{ 0,1 \}, \forall \, k \in D, \lb s_k \in S_K, s_{I(k)} \in S_{I(k)}$ and path probabilities $\pi_k(s) \in [0,1], \forall \, s \in S$.    
\end{corollary}

In particular, the objective function and constraints in Corollary \ref{cor:opt1} are linear in the decision variables $z(s_j \mid s_{I(j)})$ and the corresponding path probabilities $\pi_k(s)$. This is a MILP problem for which the optimal decision strategy can be computed with off-the-shelf MILP solvers.    

\subsection{Improving the MILP formulation}

To simplify the formulation in Corollary \ref{cor:opt1}, we note that the objective function \eqref{optimization} involves path probabilities $\pi(s)$ only for full paths $s \in S = S_{1:n}$ of length $n$. Also, the probability $\pi(s)$ of each path $s \in S$ depends on two separable components. First, for each path $s \in S$, the conditional probabilities \eqref{CondProb} of the states $s_j$ for chance nodes $j \in C$ can be multiplied to obtain the following upper bound for $\pi(s)$:
\begin{equation}\label{psn}
p(s) = \prod_{j \in C} \mathbb{P}(X_j = s_j \mid X_{I(j)} = s_{I(j)}).
\end{equation}
Second, for a given decision strategy $Z \in \mathbb{Z}$, this upper bound $p(s)$ is the actual probability of $s$ if and only if $Z$ is compatible with $s$. That is, if $z(s_j \mid s_{I(j)}) = 1$, $\forall j \in D$, the inequalities \eqref{inequality1}--\eqref{inequality4} 
imply $\pi_{j}(s) = \pi_{j-1}(s)$ for each $j \in D$. This result can be used to solve the equations \eqref{condprob}--\eqref{DecisionPi} recursively starting from $\pi_0(s) = 1$ to the last node $n$ for which $\pi_{n}(s) = p(s)$ in \eqref{psn}. Conversely, if the decision strategy $Z$ is not compatible with $s$, inequalities \eqref{inequality1}--\eqref{inequality2} imply that $\pi_n(s) \leq \pi_j(s) = 0$ for some $j \in D$. Thus, because $\pi(s) = \pi_n(s) = p(s)$ if and only if $z(s_j \mid s_{I(j}) = 1$, $\forall \, j \in D$, the optimization problem in Corollary \ref{cor:opt1} can be reformulated as
\begin{align}
\mathop{\text{maximize}}\limits_{Z \in \mathbb{Z}} ~~&  \sum_{s \in S} \pi(s) \, \mathcal{U}(s) \label{objEU} \\
\text{subject to} ~~\, & \sum_{s_j \in S_j} z(s_j \mid s_{I(j)}) = 1, &&\forall \,j \in D,  \,s_{I(j)} \in S_{I(j)} \label{eq0}\\ 
& 0 \leq \pi(s) \leq p(s),  && \forall\, s \in S \label{eq1} \\
& \pi(s) \leq z(s_j \mid s_{I(j)}), && \forall\, s\in S  \label{eq2} \\
& \pi(s) \geq p(s) + \sum_{j \in D} z(s_j \mid s_{I(j)}) - |D|,  &&\forall\, s \in S \label{eq3}\\ & z(s_j \mid s_{I(j)}) \in \{0,1\}, &&\forall\, j \in D,\, s_j \in S_j, \, s_{I(j)} \in S_{I(j)}, \label{eq4}
\end{align}
where the constraints \eqref{eq0} ensure that some decision $s_j \in S_j$ is made at each decision node $j\in D$ for every information state set $s_{I(j)} \in S_{I(j)}$ (as stated in \eqref{normconstraint}). Constraints \eqref{eq1} bound the probabilities of paths $s \in S$. Constraints \eqref{eq2} ensure that only those paths which are compatible with the decision strategy can have positive probabilities. Constraints \eqref{eq3} ensure that the probabilities of paths with negative utility $\mathcal{U}(s)$ cannot become smaller than their upper bounds $p(s)$ for paths $s$ such that $z(s_j \mid s_{I(j)}) = 1$, $j \in D$. However, because utility functions are unique to positive affine transformations and the value node has a finite number of $\prod_{i \in I(v)} |S_i|$ information states, one can choose a utility function with non-negative utilities (i.e., ${\rm min}_{s \in S} U[Y_v(s_{I(v)})] \geq  0$), which allows for omitting constraints \eqref{eq3}. For clarity, we note that in constraints \eqref{eq2}--\eqref{eq3}, the states $s_j, s_{I(j)}$ are taken from the  selected path $s \in S$. Finally, constraints \eqref{eq4} enforce the domain of all binary variables $z(s_j \mid s_{I(j)})$.

\subsection{Valid equalities}\label{sec:probcuts}

Next, we describe valid equalities to strengthen the problem formulation \eqref{objEU}--\eqref{eq4} so that it can be solved more efficiently. These equalities are derived from the problem structure and can help compute the optimal decision strategies, as shown in Section \ref{sec:CT}. However, adding these equalities directly as additional constraints may slow down the overall solution process especially for larger problems, as many of them can be derived from the problem structure. 

Alternatively, one can include these valid equalities during the solution process as ``lazy constraints'' that can be used by the MILP solver to prune nodes of the branch-and-bound tree more efficiently. One can also add them during the solution process in a cutting plane fashion as ``user cuts'' for a subset of nodes in the tree based on some criterion (or multiple criteria), for example, if the upper bound has not improved enough within some time interval. Such lazy constraints and user cuts are standard features in off-the-shelf MILP solvers.

Specifically, the first set of equalities, referred to as \emph{probability cuts}, exploit the fact that for any decision strategy $Z\in \mathbb{Z}$, the sum of the probabilities $\pi(s)$ must equal one so that $\sum\limits_{s\in S} \pi(s) = 1$. These equalities are valid for any problem that can be formulated as \eqref{objEU}--\eqref{eq4}. As an example of how a probability cut works as a lazy constraint, suppose that the optimal (fractional) solution of a node in the branch and bound tree does not satisfy the probability cut. Then, the problem at that node will be re-optimized after adding the probability cut, and if the new optimal cost is smaller than the current best primal bound, the node can be discarded. %Other cases in which the node can be pruned are if the problem becomes infeasible or if it leads to an integer solution.

The second set of equalities can be used in problems whose structure makes it possible to determine in advance for a given decision strategy $Z\in \mathbb{Z}$ how many paths $s \in S$ are \emph{active} so that  $\pi(s) > 0$. For example, if the number of active paths in any solution is $n_s$, we can define a valid equality $\sum_{s\in S} \pi(s)/p(s) = n_s$ where $p(s)$ in \eqref{psn} is the upper bound for $\pi(s)$. This approach can be generalized to asymmetric problems in which the number of active paths varies for different decision strategies. In such cases, several equalities can to be added to cover different possibilities in how the number of active paths depends on the states of decision or chance nodes. Such information, derived from an analysis of symmetries in the problem structure (see, e.g., \citealp{Bielza2011}), serve to improve computational efficiency. 

%The effectiveness of the probability cuts are demonstrated in Section \ref{sec:CT} where we provide a comparison of solution times of large problem instances with and without these cuts.
%As another example, if the number of paths in any solution belongs to the set $n_s = \{n_{s1},n_{s2},\dots,n_{sp}\}$, we can include a set of valid inequalities of the form 
%\begin{align*}
%&\sum_{s\in S} \pi(s)/p(s) = n_{si} y_i, &&\forall \, i = 1,\dots,p\\
%&y_1 + y_2 + \dots + y_s = 1 &&\\
%&y_i \in \{0,1\}, &&\forall \, i = 1,\dots,p
%\end{align*}

\section{Decision Modeling Examples}\label{sec:IE}

\subsection{Decision Programming without the No-Forgetting Assumption} \label{sec:DM}

As an example of a problem in which the \emph{no-forgetting} assumption does not hold, assume that there is an uncertain load $L$ on a built structure which can be fortified through actions $A^1$ and $A^2$ to mitigate the risk of a failure $F$ of this structure. These two decisions are informed by measurement reports $R^1$ and $R^2$ of the load $L$. The decision as to whether action $A^1$ should be implemented is informed by the report $R^1$ only and, similarly, decision $A^2$ is based on the report $R^2$ alone. In particular, the decision as to whether the fortification decision $A^1$ will be or has been installed is not known when making the decision $A^2$ (and conversely for $A^2$). The utility at the target node $T$ depends on whether or not the structure fails and how much the fortification actions cost. 

This problem structure also represents a situation where the reports are generated by sensors which inform safety controls (e.g., valves) that must activated instantaneously to prevent potential disruptions in a safety-critical system such as a nuclear plant (see, e.g., \citep{MancusoEtAl2019}). In particular, the safety must be ensured even if failures of communication equipment prevent the sensors from sharing  information with a centralised server or other sensors. 

\begin{figure}
\centering
\begin{tikzpicture}
    [decision/.style={fill=blue!80, draw, minimum size=2.5em, inner sep=2pt}, 
    chance/.style={circle, fill=orange!80, draw, minimum size=2.5em, inner sep=2pt},
    value/.style={diamond, fill=teal!60, draw, minimum size=2.5em, inner sep=2pt},
    scale=1.5]
    %\draw[step=1cm,gray,very thin] (0,0) grid (3,2);
    \node[chance]   (L) at (0, 1)  {$L$};
    \node[chance]   (La) at (0.5, 2)  {$R^1$};
    \node[chance]   (Lb) at (0.5, 0)  {$R^2$};
    \node[decision] (A) at (1.5, 2)  {$A^1$};
    \node[decision] (B) at (1.5, 0)  {$A^2$};
    \node[chance]   (F) at (2, 1)  {$F$};
    \node[value]    (T) at (3, 1)  {$T$};     
    \draw[->, thick] (L) -- (La);
    \draw[->, thick] (L) -- (Lb);
    \draw[->, thick] (L) -- (F);
    \draw[->, thick] (La) -- (A);
    \draw[->, thick] (Lb) -- (B);
    \draw[->, thick] (A) -- (F);
    \draw[->, thick] (B) -- (F);
    \draw[->, thick] (F) -- (T);
    \draw[->, thick] (A) -- (T);
    \draw[->, thick] (B) -- (T);
\end{tikzpicture}
\caption{Influence diagram of the example on double monitoring.} \label{fig:Figure1}
\end{figure}
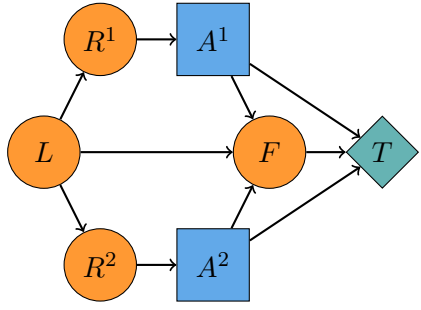

Just as in the example in Figure 12 from \citet{Zhang1994}, this problem structure is challenging in that the optimal strategy at one decision node depends on that at of the other. In particular, the no-forgetting assumption does not hold, because there is no sequence of chance nodes $C = \{ L,R^1, R^2, F \}$ and decision nodes $D = \{ A^1,A^2 \}$ such that for all decision nodes, the states of all preceding nodes would be known at the time of decision making. Figure \ref{fig:Figure1} presents an influence diagram representing this setting.

Still, this problem can be solved using \emph{Decision Programming}. The sequence $(L,R^1,R^2,A^1,A^2,F\,,T)$ captures the dependence structure: $I(R^i) = \{ L \}, \, i=1,2$ (the reports depend on the load); $I(A^i) = \{ R^i \}, \, i=1,2$ (decisions about the fortification actions are informed by respective reports); $I(F)= \{L,A^1, A^2\}$ (failure depends on the load and fortification decisions, but not on the reports); and $I(T) = \{ A^1,A^2, F \}$ (the final outcome depends on the failure and the cost of implementing the fortification actions). By using node labels to indicate sets of states for corresponding nodes, the paths are sequences $s = (l,r^1,r^2,a^1,a^2,f) \in L \times R^1 \times R^2 \times A^1 \times A^2 \times F = S$. The probabilities $p(s)$ in \eqref{psn} are $p(l,r^1,r^2,a^1,a^2,f) = \mathbb{P}(l) \mathbb{P}(r^1 \mid l) \mathbb{P}(r^2 \mid l) \mathbb{P}(f \mid  l, a^1, a^2)$, and the decision strategies are defined by $Z = (Z_1,Z_2)$ such that $Z_i : R^i \mapsto A^i$. 

Using this notation, the optimal fortification strategy can be obtained by solving the equations \eqref{eq0}--\eqref{eq4}, which in this example become
\begin{align*}
\mathop{\text{maximize}}\limits_{Z \in \mathbb{Z}} ~~&\sum_{(l,r^1,r^2,a^1,a^2,f)} \pi(l,r^1,r^2,a^1,a^2,f) U\big[Y_T(a^1,a^2,f)\big] \\
\text{subject to} ~~\, & \sum_{a^i \in A^i} z(a^i \mid r^i) = 1, && \forall~r^i \in R^i,  i=1,2 \\
& 0 \leq \pi(l,r^1,r^2,a^1,a^2,f) \leq  p(l,r^1,r^2,a^1,a^2,f), && \forall~(l,r^1,r^2,a^1,a^2,f) \in S \  \\
& \pi(l,r^1,r^2,a^1,a^2,f) \leq z(a^i \mid r^i), && \forall ~(l,r^1,r^2,a^1,a^2,f) \in S \\
& \pi(l,r^1,r^2,a^1,a^2,f) \geq p(l,r^1,r^2,a^1,a^2,f) + \sum_{i=1,2}z(a^i \mid r^i) - 2, && \forall ~(l,r^1,r^2,a^1,a^2,f) \in S \\
& z(a^i \mid r^i) \in \{0,1\}, && \forall~a^i \in A^i, r^i \in R^i, i=1,2, 
\end{align*}
where $Y_T(a^1,a^2,f)$ gives the consequences associated with the failure state $f$ and the actions $a^1$ and $a^2$. If all the decision and chance nodes have binary states, then there are altogether 8 decision variables (4 per each fortification decision) and $2^6 = 64$ paths, resulting in 4 equality constraints and 128 inequality constraints (in the second inequality constraint, the states $a^i,r^i$ are implied by the selected path and third inequality constraints can be omitted by normalizing the utility function so that it attains positive values only). 

\subsection{Decision Programming as an Extension of Contingent Portfolio Programming} \label{CPP}

\emph{Contingent Portfolio Programming} 
\citep{GustafssonSalo2005} is a methodology for determining optimal decision strategies for multi-period investment projects whose cash flows depend on (i) uncertain states of nature and (ii) project management decisions. The aim is to maximize the expected resource position at the terminal period, subject to relevant resource and consistency constraints. Risk preferences can be accounted for either by formulating risk constraints or by introducing risk measures into the objective function. CPP problems of realistic size can be tackled with off-the-shelf MILP solvers. However, a limitation of CPP is that the probabilities of the states of nature cannot depend on project decisions. Yet, such dependencies arise for instance when the projects influence market and regulatory uncertainties (for a case study, see  \citealp{Vilkkumaa2018}). 

\emph{Decision Programming} generalizes CPP models so that the probabilities of the states of nature in the scenario tree can depend endogenously on project decisions. We illustrate this by extending the example in \citet{GustafssonSalo2005} with two projects $A$ and $B$ of which one or both can be started in period 0. If a project is started, it can be continued in period 1, in which case it gives in period 2 a payoff that depends on the two chance events $C^1$ and $C^2$ which correspond to uncertain upward and downward movements in the CPP scenario tree in periods 1 and 2, respectively.

Specifically, the decision nodes $A^s, B^s$ indicate whether the projects $A$ and $B$ are started, and $A^c$ and $B^c$ correspond to decisions to continue them. The decisions to continue projects are informed by the first period movement so that $I(A^c) = \{A^s,C^1 \}, I(B^c) = \{B^s,C^1 \}$. In \citet{GustafssonSalo2005}, the second period movement depends on the first period movement (i.e., $I(C^2) = \{ C^1 \}$), but is not influenced by project decisions. 

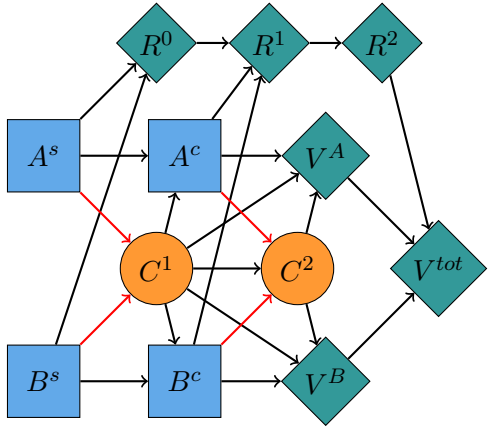
\begin{figure}
\centering 
\begin{tikzpicture}
    [decision/.style={fill=blue!80, draw, minimum size=2.5em, inner sep=2pt}, 
    chance/.style={circle, fill=orange!80, draw, minimum size=2.5em, inner sep=2pt},
    value/.style={diamond, fill=teal!80, draw, minimum size=2.5em, inner sep=2pt},
    scale=1.5]
     %\draw[step=1cm,gray,very thin] (0,0) grid (4,3);
     %\node[chance] (I) at (0, 1.5)      {$I$};
     \node[decision] (As) at (0, 2.5) {$A^s$};
     \node[decision] (Ac) at (1.25, 2.5)   {$A^c$};          
     \node[decision] (Bs) at (0, 0.5) {$B^s$};
     \node[decision] (Bc) at (1.25, 0.5)   {$B^c$};
     \node[chance]   (C1) at (1, 1.5)   {$C^1$};
     \node[chance]   (C2) at (2.25, 1.5) {$C^2$};
     \node[value]    (R0) at (1, 3.5) {$R^0$};
     \node[value]    (R1) at (2, 3.5) {$R^1$};
     \node[value]    (R2) at (3, 3.5) {$R^2$};
     \node[value]    (VA) at (2.5, 2.5) {$V^A$};
     \node[value]    (V)  at (3.5, 1.5)   {$V^{tot}$};                    
     \node[value]    (VB) at (2.5, 0.5) {$V^B$};
     %\draw[->, thick] (I) -- (As);
     %\draw[->, thick] (I) -- (Bs);
     \draw[->, thick] (As) -- (Ac);
     \draw[->, thick] (Bs) -- (Bc);
     %\draw[->, thick] (As) -- (0.5, 3.25) -- (3.5, 3.25) -- (VA);
     %\draw[->, thick] (Bs) -- (0.5, -0.25) -- (3.5, -0.25) -- (VB);
     \draw[->, thick] (As) -- (R0);
     \draw[->, thick] (Bs) -- (R0);
     \draw[->, thick] (Ac) -- (R1);
     \draw[->, thick] (Bc) -- (R1);
     \draw[->, thick] (Ac) -- (VA);
     \draw[->, thick] (Bc) -- (VB);
     \draw[->, thick] (C1) -- (Ac);
     %\draw[->, thick] (S1) -- (I);
     \draw[->, thick] (C1) -- (Bc);
     \draw[->, thick] (C1) -- (VA);
     \draw[->, thick] (C1) -- (VB);
     \draw[->, thick] (C1) -- (C2);
     \draw[->, thick, red] (As) -- (C1);
     \draw[->, thick, red] (Bs) -- (C1);
     \draw[->, thick, red] (Ac) -- (C2);
     \draw[->, thick, red] (Bc) -- (C2);
     \draw[->, thick] (VA) -- (V);
     \draw[->, thick] (VB) -- (V);
     \draw[->, thick] (C2) -- (VA);
     \draw[->, thick] (C2) -- (VB);
     \draw[->, thick] (R0) -- (R1);
     \draw[->, thick] (R1) -- (R2);
     \draw[->, thick] (R2) -- (V);
\end{tikzpicture}
\caption{The extended CPP model. The red arrows indicate additional dependencies caused by the impact of project decision on chance events $C^1$ and $C^2$.} \label{FigureExtendedCPP}
\end{figure}

Decision Programming, however, allows project decisions to impact the upward and downward movements defined by the chance events. For instance, if $C^1$ depends on $A^s$ and $B^s$ while $C^1$ depends on decisions $A^c$ and $B^c$, then $I(C^1) = \{ A^s, B^s \}$ and $I(C^2) = \{ C^1, A^c, B^c \}$. In Figure \ref{FigureExtendedCPP}, these additional dependencies are shown by the red arrows. Here, nodes $R_0,R_1,R_2$ represent resource surpluses in periods $0,1,$ and 2 after the project decisions have been taken while $V_A$ and $V_B$ represent cash flows from projects $A$ and $B$ in period 2. These cash flows are added up with $R_2$ to obtain the terminal cash position, represented by the value node $V^{\rm tot}$. 

Yet, a computational challenge with the approach of embedding project decisions in paths is that that the paths tend to become prohibitively long, because each project increases the path length by the number of decisions for the project. This challenge can be addressed by introducing decision nodes whose states represent the aggregate portfolio-level performance achieved through the selected projects. The chance events defining  the scenario tree can then be made contingent on these performance levels, while decisions concerning portfolio-level performance can be employed as constraints on the project-specific selection decisions. 

For instance, assume that the first-stage decisions specify which technology development projects will be started to generate patent-based intellectual property (P) for a platform. This intellectual property contributes subject to some uncertainties to the technical competitiveness (T) of the platform. In the second stage, it is possible to carry out application (A) development projects which, when completed, yield cash flows that depend on the market share of the platform. This market share (M) depends on the competitiveness of the platform and the number of developed applications. The aim is to maximize the cash flows from application projects less the cost of technology and application development projects.  

The structure of this problem can be modeled by introducing decision nodes $D^{P}, D^{A}$ for the development of patents and applications, respectively so that the states of these nodes correspond to continuous and contiguous intervals $d^P_i = [\underline{d}^P_i, \overline{d}^P_i), d^{A}_k = [\underline{d}^A_k, \overline{d}^A_k)$ with $\overline{d}^{P}_i =  \underline{d}^{P}_{i+1}, \overline{d}^{A}_k =  \underline{d}^{A}_{k+1} $ over ranges for the number of patents and applications that can be generated. In other words, these  states represent discretizations of these ranges, indexed with $i=1,\ldots,|D^{P}|$ and $k= 1, \ldots,|D^{A}|$.   

The technical competitiveness of the platform and its market size are represented by chance nodes $C^T$ and $C^M$ whose realizations are denoted by the states $c^T_j, c^M_l$ (where the states of these chance nodes are indicated by indices $j=1,\ldots, |C^T|$ and $l=1,\ldots,|C^M|$). The dependencies at these chance nodes can be characterized by estimating the probabilities $p(c^T_j \mid d^P_i)$ and $p(c^M_l \mid c^T_j, d^A_k)$ for the relevant combinations of states. 

At the aggregate portfolio level, the decisions consist of (i) choosing the number of patents to be generated by selecting the interval $d^P_i = [\underline{d}^P_i, \overline{d}^P_i)$ such that $z(d^P_i)=1$ and, based on this decision and the ensuing state of technical competitiveness $c^T_j$, (ii) deciding the number of applications to be developed by choosing the state $d^A_k =  [\underline{d}^A_k, \overline{d}^A_k)$ such that $z(d^A_k \mid d^P_i, c^T_j)=1$. Thus, the corresponding Decision Programming paths are sequences $s=(d^P_i, c^T_j, d^A_k, c^M_l)$ whose probabilities $\pi(s)$ are characterized by the inequalities \eqref{eq0}-\eqref{eq4} in Section 3 so that  
\begin{equation}
    \pi(s) = z(d^P_i) p(c^T_j \mid d^P_i) z(d^A_k \mid d^P_i, c^T_j) p(c^M_l \mid c^T_j, d^A_k). \label{eq:pathfordecisions}
\end{equation}

The portfolio-level decisions $d^P_i$ and $d^A_k$ can be linked to the selection of technology and application projects as follows. Assume there are $n_T$ technology projects such that the project $t \in \{ 1, \ldots, n_T \}$ requires the investment $I_t$ and produces  $O_t$ patents as its output, as well as and $n_A$ application projects such that project $a \in \{ 1, \ldots, n_A \}$ requires the investment of $I_a$ and generates $O_a$ applications. If completed, the application project $a$ provides the cash flow $V(a \mid c^M_l)$ when the market size is $c^M_l$. Using binary decision variables $x^T(t)$ and $x^A(a \mid c^T_j)$ to indicate which technology and application projects are selected (note that the competitiveness of the technology is known when starting application projects), the realized cash flow ${\cal U}(s)$ for the path $s = (d^P_i,c^T_j,d^A_k,c^M_l)$ is 
\begin{equation}
    {\cal U}(s) = \sum_{a=1}^{n_A}  x^A(a \mid d^P_i, c^T_j) \Big(V(a \mid c^M_l) - I_a \Big) - \sum_{t=1}^{n_T} x^T(t) I_t. \label{eq:cashflow}
\end{equation} 
 
The selection of technology projects is indicated by binary variables $x^T_i(t)$ which indicate whether or not the technology project $t$ is selected when the number of patents belongs to the interval $d^P_i = [\underline{d}^P_i,\overline{d}^P_i)$. Similarly, the binary variables $x^A_k(a \mid d^P_i, c^T_j)$ indicate if the application project $a$ is selected when the total number of applications is chosen to be in the interval $d^A_k = [\underline{d}^A_k, \overline{d}^A_k)$, knowing that the technological competitiveness is $c^T_j$ and that the  number of developed patents was is in the interval $d^P_i$. Because these variables indicate selections for these conditions only, we have  
\begin{eqnarray}
\sum_t x^T_i(t) & \leq & z(d^P_i) n_T, \quad \forall~i \label{eq:selectionbound1} \\
\sum_a x^A_k(a \mid d^P_i, c^T_j) & \leq &  z(d^P_i) n_A, \quad \forall~i,j,k, \label{eq:selectionbound2} \\
\sum_a x^A_k(a \mid d^P_i, c^T_j) & \leq &  z(d^A_k \mid d^P_i, c^T_j) n_A, \quad\forall~i,j,k.  \label{eq:selectionbound3}
\end{eqnarray} 

To ensure that the selection of the technology and application projects are aligned with decisions $z(d^P_i)$ and $z(d^A_k \mid  d^P_i, c^T_j)$, we have the  constraints   
\begin{eqnarray}
  & \underline{d}^P_i - \big[1 - z(d^P_i) \big] M  \leq \sum_{t} x^T_i(t) O_t  \leq  \overline{d}^{P}_i + [1 - z(d^P_i)]M - \varepsilon  \label{eq:platform1} \\
 & \underline{d}^A_k - \big[1 - z(d^A_k \mid  d^P_i, c^T_j) \big] M  
\leq \sum_{a} x^A_k(a \mid  d^P_i, c^T_j) O_a \leq \overline{d}^A_k + [1 - z(d^A_k \mid  d^P_i, c^T_j)]M - \varepsilon, 
\label{eq:platform2} 
\end{eqnarray}
where $M$ is a large constant and the terms for $\varepsilon = \frac12 {\rm min}_{t, a} \{ O_t,O_a \}$ eliminate the possibility of making selections which lead to the upper bound of the interval. By constraint \eqref{eq:platform1}, the number of patents is in the selected interval $d^P_i = [\underline{d}^P_i, \overline{d}^P_i)$ for which $z(d^P_i) = 1$, while  \eqref{eq:platform2} ensures that the number of applications is in the interval $d^A_k = [ \underline{d}^A_k, \overline{d}^A_k)$ when $z(d^A_k \mid d^P_i, c^T_j) = 1$.  

By construction, the terms in the cash flow expression \eqref{eq:cashflow} can be written as $x^T(t) = \sum_{i} x^T_i(t) z(d^P_i)$  and $x^A(a \mid d^P_i, c^T_j) = \sum_{k} x^A_k(a \mid d^P_i, c^T_j) z(d^A_k \mid d^P_i, c^T_j)$. Using these and noting \eqref{eq:pathfordecisions}, each term $\pi(s) {\cal U}(s)$ in \eqref{optimization} can be written as  
\begin{eqnarray*}
&& \sum_{i,j,k,l} z(d^P_i) p(c^T_j \mid d^P_i) z(d^A_k \mid d^P_i, c^T_j) p(c^M_l \mid c^T_j, d^A_k) \\
&& \Big[ \sum_{a}  \Big\{\sum_{k'} x^A_{k'}(a| d^P_i, c^T_j) z(d^A_{k'}|d^P_i, c^T_j) \Big\} \Big(V(a \mid c^M_l) - I_a \Big) - \sum_{t} \Big\{ \sum_{i'} x^T_{i'}(t) z(d^P_{i'}) \Big\} I_t\Big] \\
&& = \sum_{i,j,k,l} z(d^P_i) p(c^T_j \mid d^P_i) z(d^A_k \mid  d^P_i, c^T_j) p(c^M_l \mid c^T_j, d^A_k) 
\Bigg[ \sum_{a} x^A_k(a| d^P_i, c^T_j) \Big(V(a \mid c^M_l) - I_a \Big) - \sum_{t} x^T_i(t) I_t\Bigg] \\
&& = \sum_{i} \Bigg\{ \sum_{j,k,l}  p(c^T_j \mid d^P_i)  p(c^M_l \mid c^T_j, d^A_k) 
\Bigg[ \sum_{a} x^A_k(a|d^P_i, c^T_j)  \Big(V(a \mid c^M_l) - I_a \Big) \Bigg] - \sum_{t} x^T_i(t) I_t \Bigg\}, \label{eq:platformobjective}
\end{eqnarray*}
because $z(d^P_i) z(d^P_{i'}) = 0, i \neq i'$ and $z(d^A_k  \mid d^P_i, c^T_j) z(d^A_{k'} \mid d^P_i, c^T_j) = 0, k \neq k'$. By constraints \eqref{eq:selectionbound1}-\eqref{eq:selectionbound3}, the project selections $ x^T_i(t)$ and $x^A_k(a|d^P_i, c^T_j)$ can affect the objective function only for the selected intervals $z(d^P_i) = z(d^A_k \mid d^P_i, c^T_j) = 1$ while the constraints \eqref{eq:platform1}-\eqref{eq:platform2} imply that the portfolios selections are indeed aligned with these selections. Thus, we have an MILP problem involving both the  decision variables $z(d^P_i)$ and $z(d^A_k \mid d^P_i, ^T_j)$, as well as the selection of technology and application projects $x^T_i(t), x^A_k(a \mid d^P_i, c^T_j)$.

Importantly, in this problem structure the number of paths $|D^P| \times |C^T| \times |D^A| \times |C^M|$  stays the same regardless of the size of the sets from which technology and application development projects are selected. Nor does the number of constraints grow with number of project candidates. This helps solve much bigger problems while still accounting for the endogenous impact that the selected projects have on technological competitiveness and market share. In effect, optimization becomes indispensable when there are more candidate projects to be selected from. For instance, from a set of $30$ candidate projects one can build $2^{30} \approx 1,074 \times 10^9$ portfolios, wherefore explicit enumeration is no longer a viable approach. A further benefit of this layered structure is that the elicitation of conditional probabilities for the chance nodes can be largely separated from the consideration of individual technology and application projects. Still, the compatibility of the two layers is ensured by the optimization model which can readily accommodate many kinds of constraints, such as those representing budgetary constraints, risk preferences or logical dependencies.

%The above reformulation is more expressive than the standard CPP formulation, which does not allow the probabilities of the second period movement $\mathbb{P}(c^2 \mid c^1)$ to depend on project decisions $A^c,B^c$, as represented in Figure \ref{FigureExtendedCPP}. In contrast, by Theorem \ref{scenariopaths} and \eqref{condprob} we have $\mathbb{P}(a^s, b^s, c^1, a^c, b^c, c^2|Z) = \pi(a^s, b^s, c^1, a^c, b^c, c^2) = \mathbb{P}(c^2 \mid c^1, a^c, b^c) \pi_5(a^s, b^s, c^1, a^c, b^c)$ in which the conditional probability in the first term allows such dependencies to be modeled. Such dependencies would occur, for instance, when developing applications for a technology platform whose market penetration depends on the cumulative level of investments into these applications. In this case, this market size could be modeled as a chance node with a range of discrete states so that the probabilities of states representing greater penetration levels would become higher with increased investments.  

\section{Extensions to modeling chance constraints and multiple value nodes}

We next extend the formulation in Section 3 through approaches for modeling risk measures and chance constraints, as well as for computing non-dominated strategies in problems which have multiple objectives represented by different value nodes.

Apart from the use of nonlinear utility functions $U(\,\cdot\,)$ in \eqref{preutil}, risk preferences can be accounted through risk measures $\rho$ that map decision strategies to non-negative real numbers and can be introduced as additional terms into the objective function or employed as constraints. In the following, we assume that, at the value nodes $v \in V$, the aim is to maximize the consequences $\mathcal{C}(s) = Y_v(s_{I(v)}) \in \mathbb{C}$, which are assessed using real numbers.  

\subsection{Absolute and lower-semi absolute deviation}

Let $t \in \mathbb{R}$ be a given target level for consequences and define the non-negative deviation variables 
\begin{equation}
\label{DeltaT}
\Delta^+_t(s) = \max\{0,\mathcal{C}(s) - t \}, \qquad 
\Delta^-_t(s) = \max\{0,t - \mathcal{C}(s) \}.
\end{equation}
By construction, $\Delta^+_t(s)$ (respectively $\Delta^-_t(s)$) measures how much the consequence $\mathcal{C}(s)$ is above (below) the target level $t$. The deviations \eqref{DeltaT} can be precomputed for the information states $S_{I(v)}$ at the value node $v$. The expected downside risk (EDR) of a decision strategy $Z \in \mathbb{Z}$ relative to the target level $t$ is 
\begin{equation}\label{eq:EDR}
\rho_{\textrm{EDR}}(Z;t) = \sum_{s \in S} \pi(s) \Delta^-_t(s).
\end{equation}
If $t$ is chosen to be the expected value of consequences $\mathbb{E}[\mathcal{C} \mid Z] = \sum_{s\in S} \pi(s) Y_v(s_{I(v)})$ for the decision strategy $Z$, the corresponding non-negative deviation (decision) variables $\Delta^+_{\mathbb{E}[\mathcal{C} \mid Z]}(s),  \Delta^-_{\mathbb{E}[\mathcal{C} \mid Z]}(s)$ can be employed with the constraint 
$$\mathcal{C}(s)  - \Delta^+_{\mathbb{E}[\mathcal{C} \mid Z]}(s) + \Delta^-_{\mathbb{E}[\mathcal{C} \mid Z]}(s) =  
\mathbb{E}[\mathcal{C} \mid Z]$$ 
to capture the deviations from $\mathbb{E}[\mathcal{C}]$. The absolute deviation (AD) and the lower semi-absolute deviation (LSAD) are then given by 
\begin{eqnarray}
\label{eq:AD}
\rho_{\textrm{AD}}(Z) & = & \sum_{s \in S} \pi(s) \big[\Delta^+_{\mathbb{E}[\mathcal{C} \mid Z]}(s) + \Delta^-_{\mathbb{E}[\mathcal{C} \mid Z]}(s) \big] \\
\label{eq:LSAD} 
\rho_{\textrm{LSAD}}(Z) & = & \sum_{s \in S} \pi(s) \Delta^-_{\mathbb{E}[\mathcal{C} \mid Z]}(s).
\end{eqnarray}

These measures can be used to augment the objective function through an additional additive term which penalizes for risk. For example, if the aim is to maximize expected consequences while accounting for risks through (lower semi-)absolute deviation, one possibility is to formulate the objective function as $\max_{Z \in \mathbb{Z}} \big{\lbrace} (1-\phi)\mathbb{E}[X_v|Z] - \phi \rho_{\textrm{LSAD}}(Z) \big{\rbrace}$ where $\phi \in [0,1]$ is a weighting coefficient that reflects the DM's risk aversion. Alternatively, as an example of using risk measures to constrain feasible decision strategies, assume that the consequences are defined as profits reported in kUSD. Then the constraint $\rho_{\textrm{AD}}(Z) \leq 10$ would rule out any strategy $Z \in \mathbb{Z}$ for which profits can be expected to differ more than 10 kUSD from the level of its expected profits $\mathbb{E}[\mathcal{C} \mid Z]$. 

\subsection{Chance constraints and Value-at-Risk}

Probabilistic chance constraints can be modeled as linear inequalities on the path probabilities $\pi(s)$ which depend linearly on the  decision variables. For example, to assess whether the consequences $\mathcal{C}(s)$ meet or exceed the stated target level $t \in \mathbb{R}$, we define the parameters
\begin{equation}\label{XTDefinition}
\Lambda_t(s)= \begin{cases}
1, & \quad \textrm{if~} \mathcal{C}(s) \geq t \\
0, & \quad \textrm{otherwise}.
\end{cases} 
\end{equation}
If the outcome is required to reach the target level $t$ with a probability that is higher than or equal to a stated threshold level $p_t$, we can impose the constraint   
\begin{equation} \label{ChanceConstraint}
\mathbb{P}\Big( \{ s \mid \mathcal{C}(s) \geq t \}  \mid Z \Big) = \sum_{s \in S} \pi(s) \Lambda_t(s) \geq p_t,
\end{equation}
which is linear in the path probabilities $\pi(s)$. The terms $\Lambda_t(s), \forall \, s \in S$, need to be defined only for $s_{I(v)} \in S_{I(v)}$. 

In the present context where the probability distributions over consequences are discrete, the Value-at-Risk (VaR) risk measure for the strategy $Z$ can be defined as 
\begin{equation} \label{VaR}
{\rm VaR}_{\alpha}(Z) =  F^{-1}_Z(\alpha) = \sup \braces{ t \mid \mathbb{P}(s \mid \mathcal{C}(s) \leq t) < \alpha},
\end{equation}
where $F^{-1}_Z$ is the inverse function of the cumulative probability distribution $F_Z : \mathbb{C} \mapsto [0,1]$ which is defined as $F_Z(t) = \sum_{s \mid \mathcal{C}(s) \leq t} \pi(s)$. 

Because the probability distribution over the set of paths is discrete, the definition \eqref{VaR} means that consequences which are less than or equal to ${\rm VaR}_\alpha(Z)$ can occur with a probability greater than $\alpha$ \citep{RockafellarUryasev2002}. This is the case if ${\rm VaR}_{\alpha}(Z)$ coincides with a consequence where the cumulative probability distribution function jumps from a level below $\alpha$ to one that exceeds $\alpha$ so that $\mathbb{P}\Big(\{ s \mid \mathcal{C}(s) < {\rm VaR}_\alpha(Z)\} \Big) < \alpha < \mathbb{P}\Big(\{s \mid \mathcal{C}(s) \leq {\rm VaR}_{\alpha}(Z)\} \Big)$.  

Constraints such as \eqref{ChanceConstraint} can be employed to introduce VaR requirements. That is, if the probability $\alpha > 0$ is associated with the corresponding VaR level $t^\alpha_{\rm VaR}$, then the path probabilities of feasible decision strategies must satisfy the constraint
\begin{equation} \label{VaRConstraint}
\sum_{s \in S} \pi(s) \Lambda_{t^\alpha_{\rm VaR}}(s) \le \alpha. 
\end{equation}
This approach can be generalized to introduce chance constraints on the states of nodes $k \in C \cup D$ as well. For instance, assume that the state at node $k$ needs to be in some set $\tilde S_k \subset S_k$ with a probability which is less than or equal to $\tilde p_k$. This requirement can be represented by the constraint $\sum_{s \in S} \pi(s) \Lambda_{\tilde S_k}(s) \geq \tilde p_k$ where $\Lambda_{\tilde S_k}(s) = 1$ if $s_k \in \tilde S_k$ and $\Lambda_{\tilde S_k}(s) = 0$ otherwise. Thus, for example, for a decision node $k \in D$, one could require that the probability of having to employ extraordinary decisions, as represented by the states $\tilde S_k$, does not exceed a pre-specified probability $\tilde p_k$. 

\subsection{Conditional Value-at-Risk}

For a given probability $\alpha > 0$ and a decision strategy $Z\in \mathbb{Z}$, the Conditional Value-at-Risk (CVaR) is the expected level of consequences in the event that the realized consequence is in the $\alpha \in (0,1]$ lower tail of the probability distribution. Contributions to this expectation come from (i) paths $s \in S^<_{{\rm VaR}_\alpha(Z)} = \{ s \in S ~|~ \mathcal{C}(s) < {\rm VaR}_\alpha(Z) \}$ which lead to consequences strictly less than ${\rm VaR}_\alpha(Z)$ and (ii) paths $s \in S^=_{{\rm VaR}_\alpha(Z)} = \{ s \in S ~|~ \mathcal{C}(s) = {\rm VaR}_\alpha(Z) \}$ which lead to the consequence ${\rm VaR}_\alpha(Z)$. The share of the probability of these latter paths that needs to be accounted in the computation of the CVaR level is the  difference $\alpha - \mathbb{P}( \{s ~|~ \mathcal{C}(s) < {\rm VaR}_\alpha(Z)\}) = \alpha - \sum_{s \in S^<_{{\rm VaR}_\alpha(Z)}}{\pi(s)}$. Thus, as in \citet{LiesioSalo2012}, we define the risk measure $\textrm{CVaR}_{\alpha}(Z)$ as  
\begin{equation} \label{CVaR}
 \textrm{CVaR}_{\alpha}(Z) = \frac1{\alpha} \left(\sum_{s \in S^{<}_{{\rm VaR}_\alpha(Z)}} \pi(s)\mathcal{C}(s) + \sum_{s \in S^=_{{\rm VaR}_\alpha(Z)}} \left(\alpha - \sum_{s \in S^<_{{\rm VaR}_\alpha(Z)}} \pi(s)\right)\mathcal{C}(s)\right).
\end{equation}

By Proposition \ref{VarProposition}, the VaR and CVaR levels for a given probability level $\alpha > 0$ and decision strategy $Z \in \mathbb{Z}$ can be determined by solving the optimization problem \eqref{beginCVaR}--\eqref{endCVaR} with precomputed parameters $c^* = \max \{ \mathcal{C}(s) \mid s \in S \}, c^\circ =  \min \{ \mathcal{C}(s) \mid s \in S \}, M = c^*-c^\circ$ and $\epsilon = \frac12 \min \{ |\mathcal{C}(s) - \mathcal{C}(s')| \mid  |\mathcal{C}(s) - \mathcal{C}(s')| > 0, s, s' \in S \}$. 

\begin{proposition} \label{VarProposition}
Choose $\alpha \in (0,1]$ and let $\pi(s), \forall \, s \in S$, be the path probabilities for a decision strategy $Z \in \mathbb{Z}$. Then the optimization problem  
\begin{align}
& {\rm min} \qquad \eta \\
& \eta - \mathcal{C}(s) \leq M\lambda(s), && \forall s \in S \label{beginCVaR} \\
& \eta - \mathcal{C}(s) \geq (M + \epsilon)\lambda(s) - M, && \forall s \in S \label{CVaR2} \\
& \eta - \mathcal{C}(s) \leq (M + \epsilon)\overline{\lambda}(s) - \epsilon, && \forall s \in S \label{CVaR3} \\
& \eta - \mathcal{C}(s) \geq M(\overline{\lambda}(s) - 1), && \forall s \in S \label{CVaR4} \\
& \overline{\rho}(s) \leq \overline{\lambda}(s), && \forall s \in S  
\label{CVaR5} \\
& \pi(s) - (1-\lambda(s)) \leq \rho(s) \leq \lambda(s), && \forall s \in S 
\label{CVaR6}\\
& \rho(s) \leq \overline{\rho}(s) \leq \pi(s), && \forall s \in S 
\label{CVaR7} \\
& \sum_{s \in S} \overline{\rho}(s) = \alpha, \label{CVaR8} \\
& \overline{\lambda}(s), \lambda(s) \in \{0,1\}, &&\forall\, s\in S \\
& \overline{\rho}(s), \rho(s) \in [0,1], &&\forall\, s\in S    \\
& \eta \in [c^\circ, c^*], \label{endCVaR} 
\end{align}
has a solution such that the optimum value $\eta^* = {\rm VaR}_{\alpha}(Z)$ and ${\rm CVaR}_{\alpha}(Z) = \dfrac1 {\alpha}{\sum_{s \in S} \overline{\rho}(s)\mathcal{C}(s)} $. 
\end{proposition}

{\noindent \sc Proof}. See Appendix A. 

An inspection of the proof of Proposition \ref{VarProposition} shows that for any feasible solution to the constraints \eqref{beginCVaR}--\eqref{endCVaR}, the expression $\sum_{s \in S} \overline{\rho}(s) \mathcal{C}(s)/ \alpha$ gives the correct ${\rm CVaR}_{\alpha}(Z)$ risk measure for $Z$. Thus, if the expectation of consequences in the lower $\alpha$-tail of the probability distribution over consequences is required to be greater than or equal to the lower bound $t^{\alpha}_{\rm CVaR}$, this requirement can be enforced by adding the constraints \eqref{beginCVaR}-\eqref{endCVaR} and $\sum_{s \in S}\overline{\rho}(s)\mathcal{C}(s) \geq \alpha t^\alpha_{\rm CVaR}$ to \eqref{eq0}-\eqref{eq4}.  

One approach to address trade-offs between the maximization of conditional expectations for different levels of $\alpha$ is to treat these as different objectives with respective weighting coefficients. Thus, combining the unconditional expectation with the selected $\alpha \in (0,1)$ for CVaR leads to the problem 
\begin{eqnarray}
\mathop{\text{maximize}} \limits_{Z \in \mathbb{Z}} ~~& w  \bigg( \sum_{s \in S} \pi(s) \mathcal{C}(s)  \bigg) + (1-w)\bigg( \frac1{\alpha}\sum_{s \in S}\overline{\rho}(s)\mathcal{C}(s) \bigg) \label{ObjFunCVaR}\\
\text{subject to} ~~& \eqref{eq0}-\eqref{eq4}, \eqref{beginCVaR}-\eqref{endCVaR} 
\end{eqnarray}
where the parameter $w \in (0,1)$ represents trade-offs between (i) the overall expectation in the first term of \eqref{ObjFunCVaR} and (ii) the expectation in the lower $\alpha$-tail as expressed by the second term. Specifically, the ratio $\dfrac{1 - w}{w}$ indicates how much of the overall expectation the DM is willing to give up in return for improving the CVaR level by one unit, regardless of the overall expectation. 

\subsection{Multiple value nodes and objectives} \label{sec:multiple}

The consideration of CVaR levels together with the maximization of expected consequences is an example of the more general case with multiple objectives $n_V > 1$. In this case, attention can be focused on non-dominated strategies $Z \in \mathbb{Z}_{ND}$ such that there is no other feasible strategy $Z' \in \mathbb{Z}_{F}$ whose expectation is equal to or higher than that of $Z$ at each value node and strictly higher for at least one value node, i.e.,  
$$
Z \in \mathbb{Z}_{ND} \iff Z \in \mathbb{Z}_F \land \not\exists Z' \in \mathbb{Z}_F {\rm~such~that~}  \mathbb{E}[\mathcal{C}_v \mid Z'] \geq \mathbb{E}[\mathcal{C}_v \mid Z], \ \forall v \in V,
$$
where $\mathbb{E}[\mathcal{C}_v \mid Z] = \sum_{s \in S} \pi(s)\mathcal{C}_v(s)$ denotes the expectation at value node $v \in V$ and the inequality is strict for at least one value node $v \in V$. 

Because the strategies are choices from a discrete set of alternatives, this is a discrete multi-objective optimization problem (MOO) in which the objectives correspond to the maximization of expectations for different value nodes. Thus, it can be solved with algorithms for this problem class. \citet{HolzmannSmith2018} provide an extensive review and propose an algorithm based on augmented Tchebychev norm, in which choices about the initial step size need to be made. 

The weighting approach in \eqref{ObjFunCVaR} or, more generally, the maximization of the expression $\sum_{v \in V} w_v \mathbb{E}[\mathcal{C}_v \mid Z]$ can be employed to generate non-dominated strategies. However, a shortcoming of the weighting approach is that it does not necessarily generate all non-dominated strategies even if all non-negative weighting coefficients $w_v \geq 0, \forall \, v \in V$ such that $\sum_{v \in V} w_v = 1$, are employed. This will be the case if a non-dominated strategy $Z' \in Z_{ND}$ is dominated by a weighted linear combination of other non-dominated strategies $Z_1,\ldots, Z_k \in Z_{ND}$ so that for some selection of positive weights $\omega_i > 0$ with $\sum_{i=1}^k \omega_i = 1$, it holds that $\mathbb{E}[\mathcal{C}_v \mid Z'] \leq \sum_{i=1}^k \omega_i \mathbb{E}[\mathcal{C}_v \mid Z_i]$ for all $v \in V$ (with a strict inequality for some $v \in V$). 

This notwithstanding, the weighting approach can be adapted to generate all non-dominated strategies. First, if $Z' \in \mathbb{Z}_{ND}$ is a non-dominated strategy, then it can be excluded when computing further candidates for non-dominated strategies by introducing the linear constraint
\begin{equation} \label{z_exclusion}
\sum_{\{ (s_i,s_{I(i)}) \mid z'(s_i \mid s_{I(i)}) = 0 \}} z(s_i \mid s_{I(i)}) + \sum_{d \in D} \prod_{i \in I(d)} |S_i| -  \sum_{\{ (s_i,s_{I(i)}) \mid z'(s_i \mid s_{I(i)}) = 1 \}}  z(s_i \mid s_{I(i)}) \geq 1, 
\end{equation}
where $z'(s_i \mid s_{I(i)}), s_i \in S_i, s_{I(i)} \in S_{I(i)}$ are the decision variables for $Z' \in \mathbb{Z}^*$. In \eqref{z_exclusion}, the left side for strategy $Z$ will be greater than one if and only if $Z$ differs from $Z'$. 

Second, if $Z' \in Z_{ND}$, then further candidates for non-dominated strategies must not be dominated by $Z'$. A necessary condition for this can be formulated by defining the binary variables $\lambda_{Z',v}^+(Z),$ $\lambda_{Z',v}^-(Z) \in \{0,1 \}, \forall \, v \in V$ so that $\lambda_{Z',v}^+(Z) + \lambda_{Z',v}^-(Z) = 1 $ and
\begin{eqnarray}
\mathbb{E}[U_v \mid Z] & \leq & \mathbb{E}[U_v \mid Z'] +  M \lambda_{Z',v}^+(Z) \label{lambdaZv+} \\
\mathbb{E}[U_v \mid Z'] & \leq & \mathbb{E}[U_v \mid Z] +  M \lambda_{Z',v}^-(Z) \label{lambdaZv-} 
\end{eqnarray}
where $M$ is a large constant (e.g., slightly greater than $c^* = {\rm max}_{s \in S} \mathcal{C}(s)$). Now, consider any solution to \eqref{lambdaZv+}--\eqref{lambdaZv-} such that $\lambda_{Z',v}^+(Z) = 0, \forall \, v \in V$. Then $\mathbb{E}[U_v \mid Z]$ is either strictly less than $\mathbb{E}[\mathcal{C}_v \mid Z]$ for all $v \in V$ so that $Z$ is dominated by $Z'$; or if not, there exists some $v' \in V$ such that $\mathbb{E}[\mathcal{C}_{v'} \mid Z] = \mathbb{E}[\mathcal{C}_{v'} \mid Z']$ so that the values of the variables $\lambda_{Z',{v'}}^-(Z) = 1, \lambda_{Z',{v'}}^+(Z) = 0$ can be switched to $\lambda_{Z',{v'}}^-(Z) = 0, \lambda_{Z',{v'}}^+(Z) = 1$, in which case the constraints \eqref{lambdaZv+}-\eqref{lambdaZv-} are still satisfied. Thus, for any strategy $Z$ which is not dominated by $Z'$ there will exist a solution such that 
\begin{equation}
\sum_{v \in V} \lambda_{Z',v}^+(Z) \geq 1, \ Z' \in Z_{ND}. \label{lambda+_sum}
\end{equation}

The above constraints \eqref{lambdaZv+}--\eqref{lambdaZv-} and \eqref{lambda+_sum} constitute a necessary but not a sufficient condition. That is, it is possible that the candidate solution $Z^{''}$ which maximizes $\sum_{v \in V} w_v \mathbb{E}[\mathcal{C}_v \mid Z]$ is dominated by $Z'$, because it is possible that the value nodes can be partitioned into non-empty sets $V^{=} \cup V^{<} = V $ such that $\mathbb{E}[\mathcal{C}_v \mid Z^{''}] = \mathbb{E}[\mathcal{C}_v \mid Z'], v \in V^{=}$ and $\mathbb{E}[\mathcal{C}_v \mid Z^{''}] < \mathbb{E}[\mathcal{C}_v \mid Z'], v \in V^{<}$, i.e., $Z^{''}$ is dominated by $Z'$. Consequently, explicit dominance checks are needed to evaluate whether the candidate solution $Z^{''}$ is non-dominated. If it is, the set of non-dominated strategies can be updated by adding $Z^{''}$ to this set and by introducing the constraint \eqref{z_exclusion} to eliminate $Z^{''}$ from further consideration. Adding this constraint to \eqref{lambdaZv+}--\eqref{lambdaZv-} for $Z^{''}$ does not prevent the computation of alternative strategies whose expectations are the same for all value nodes, as such strategies do not dominate each other. 

Next, the algorithm can be iterated by maximizing $\sum_{v \in V} w_v \mathbb{E}[\mathcal{U}_v \mid Z]$ to generate further candidate strategies and augmenting the sets of non-dominated strategies and constraints accordingly. Because the number of non-dominated strategies is finite, the algorithm will generate them all by construction. In general, the algorithm is likely to perform well if the generation of non-dominated strategies in the early stages helps eliminate many non-dominated strategies. This would be the case, for instance, if there is a strong positive correlation between the expectations of different value nodes.   

\section{Computational experiments}\label{sec:CT}

We next report results from computational experiments to demonstrate the practical viability of Decision Programming and illustrate how its performance scales with increasing problem size. All implementation were coded in Julia 1.1.0, using the package JuMP to implement models which were solved with Gurobi 8.1.0.

\subsection{N-monitoring instances}

The $N$-monitoring problem has the same structure as the double monitoring problem in Section \ref{sec:DM} except that there are $N$ binary reinforcement decisions of which each is informed by its own load report with possible states {\it low} and {\it high}. For every problem size, we solve 100 instances with randomly generated data, both with and without the probability cuts described in Section \ref{sec:probcuts}. 

Data sets with plausible characteristics were generated as follows. The utility of the structure not failing was set to $100$ and that of failing to $0$. For the load node $L$, the probability of the high load state was generated from the uniform distribution $U(0,1)$ over the unit interval and the remaining probability was assigned to the low load state. For each load level and report, the probability of receiving a correct report was taken to be ${\rm max} \{ x, 1- x\}$ where $x$ was generated from the uniform distribution $U(0,1)$. Further realizations of $x,y$ from $U(0,1)$ were used to set the prior probability of failure in the case of high load to ${\rm max} \{ x, 1- x\}$ and that in the case of low load to ${\rm min} \{ y, 1- y\}$. The costs of fortification $c_i, i=1,\ldots,N$ actions were also generated from $U(0,1)$. The posterior probability of failure after implementing the actions $A \subseteq \{ 1, \ldots, N \}$ was taken to be that of the prior divided by $e^{\sum_{i \in A} c_i}$, meaning that these actions could only decrease the probabilities of failure and that the more costly actions would be more effective in doing so. In particular, this is an example of a portfolio problem with endogenous uncertainties in which the probability of failure is impacted by the {\it portfolio} of fortification decisions.   

\begin{table}[htp!] %{tab:n-m_solution_time}
	\centering
	\caption{Results on the 100 randomly generated $N$-monitoring instances.}\label{tab:n-m_solution_time}
	\setlength{\tabcolsep}{5.0pt}
	{\footnotesize{
\begin{tabular}{c@{\extracolsep{35pt}}c@{\extracolsep{0pt}}r@{\extracolsep{35pt}}r@{\extracolsep{0pt}}r@{\extracolsep{35pt}}r@{\extracolsep{0pt}}r}
\toprule
\multicolumn{1}{l}{} & 
\multicolumn{2}{c}{Number of variables} &
\multicolumn{2}{c}{No probability cuts} &
\multicolumn{2}{c}{With probability cuts} \\
\cmidrule{2-3}
\cmidrule{4-5}
\cmidrule{6-7}
\multicolumn{1}{l}{\# Nodes} &
\multicolumn{1}{r}{Binary} &
\multicolumn{1}{r}{Real} &
\multicolumn{1}{r}{A} & 
\multicolumn{1}{r}{SD} & 
\multicolumn{1}{r}{A} & 
\multicolumn{1}{r}{SD}\\
\midrule
2  &  \p{1}8  &   64~~  &	    0.01 	&	 0.01   &	   0.01 	&	  0.00 \\
3  & 12  &     256~~  &	    0.12 	&	 0.08 	&	   0.02 	&	  0.01 \\
4  & 16  &    1~024~~  &	    0.79 	&    0.53	&      0.07 	&	  0.02 \\
5  & 20  &    4~096~~  &     5.94 	&    2.80   &      0.35  	&	  0.19 \\
6  & 24  &   16~384~~  &    77.35    &   46.31  	&	   2.44 	&	  1.63 \\
7  & 28  &   65~536~~  &   676.35    &  468.09  	&     20.58 	&    17.48 \\
8  & 32  &  262~144~~  &	8~474.00 & 7~377.28  &	 268.93 	&	330.89 \\
9  & 36  & 1~048~576~~  & 	  -      &	 	-    &   1~727.19    &  2~880.20 \\
%10	&	 	&	 	&	 	&	 		\\
\bottomrule
\end{tabular}
}}
\end{table}

Table \ref{tab:n-m_solution_time} shows the computational times in seconds needed to solve the randomly generated instances, comparing the computational performance with and without the probability cuts discussed in Section \ref{sec:probcuts}. The results are provided in terms of the average (A) and standard deviation (SD) among 100 replications. A time limit of 25~200 seconds (7 hours) was imposed to all experiments. The entry ``-'' denotes cases for which no solution could be found within the 7h time limit. In the case with $N=9$, for example, there are $4^9 = 262~122$ different decision strategies. As can be observed, the inclusion of the probability cuts greatly improve the performance of the solver.

\subsection{The pig farm problem}

\begin{figure}[h]
\centering 
\begin{tikzpicture}
    [decision/.style={fill=blue!80, draw, minimum size=2.5em, inner sep=2pt}, 
    chance/.style={circle, fill=orange!80, draw, minimum size=2.5em, inner sep=2pt},
    value/.style={diamond, fill=teal!60, draw, minimum size=2.5em, inner sep=2pt},
    scale=1.5]
     %\draw[step=1cm,gray,very thin] (0,0) grid (4,3);
     \node[chance] (h1) at (0, 3)      {$h_1$};
     \node[chance] (h2) at (1, 3)      {$h_2$};
     \node[chance] (h3) at (2, 3)      {$h_3$};
     \node[chance] (h4) at (3, 3)      {$h_4$};
     \node[value]  (u4) at (4, 3)      {$u_4$};
     \node[chance] (t1) at (0, 2)      {$t_1$};
     \node[chance] (t2) at (1, 2)      {$t_2$};
     \node[chance] (t3) at (2, 2)      {$t_3$};
     \node[decision] (d1) at (0, 1)    {$d_1$};
     \node[decision] (d2) at (1, 1)    {$d_2$};
     \node[decision] (d3) at (2, 1)    {$d_3$};
     \node[value] (u1) at (0, 0)       {$u_1$};
     \node[value] (u2) at (1, 0)       {$u_2$};
     \node[value] (u3) at (2, 0)       {$u_3$};
     \draw[->, thick] (h1) -- (t1);
     \draw[->, thick] (h1) -- (h2);
     \draw[->, thick] (h2) -- (t2);
     \draw[->, thick] (h2) -- (h3);
     \draw[->, thick] (h3) -- (t3);
     \draw[->, thick] (h3) -- (h4);
     \draw[->, thick] (h4) -- (u4);
     \draw[->, thick] (t1) -- (d1);
     %\draw[->, thick] (t1) -- (d2);
     %\draw[->, thick] (t1) -- (d3);
     \draw[->, thick] (t2) -- (d2);
     %\draw[->, thick] (t2) -- (d3);
     \draw[->, thick] (t3) -- (d3);
     %\draw[->, thick] (d1) -- (d2);
     \draw[->, thick] (d1) -- (h2);
     \draw[->, thick] (d1) -- (u1);
     \draw[->, thick] (d2) -- (h3);
     \draw[->, thick] (d2) -- (u2);
     %\draw[->, thick] (d2) -- (d3);
     %\draw[->, thick] (d1) -- (0.5, 0.5) -- (1.5,0.5) -- (d3);
     \draw[->, thick] (d3) -- (h4);
     \draw[->, thick] (d3) -- (u3);
\end{tikzpicture}
\caption{The pig farm problem with three decision periods \citep{LauritzenNilsson2001}.} \label{FigurePigs}
\end{figure}
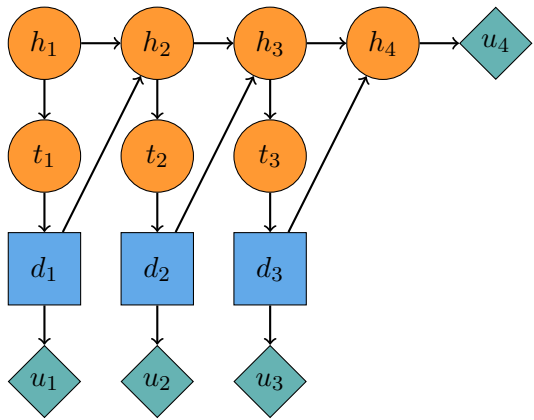

In the pig farm problem (for details, see \citealp{LauritzenNilsson2001}), a veterinary doctor visits a pig farm each month to test each pig for a disease and decides, based on the uncertain test result, whether or not to inject the pig with a drug which has both curing and preventive effects but which comes at a cost. After four months, the pigs are sold, whereby healthy pigs command a higher market price than diseased ones. There is no access to individual records for each pig, and thus the doctor has to make the treatment decision based on the age of the pig and the most recent test result. This problem is represented by a limited memory influence diagram (LIMID) in Figure \ref{FigurePigs}.

Despite its practical relevance and conceptual simplicity, this problem is not \emph{soluble} (meaning that the global optimum it not obtainable by means of solving a sequence of local optimization problems formulated for the different nodes; for details, see \citealp{LauritzenNilsson2001}, Definition 14) and consequently the Single Policy Update method proposed by the authors is not guaranteed to find globally optimal solutions. With Decision Programming, this problem can nevertheless be modeled and solved to global optimality rather efficiently. Table \ref{tab:pigs_solutions} presents the optimal solutions and their computation times both for the original four-month version of the problem with three decision periods (in which there are 64 different strategies, corresponding to $4 \times 4 \times 4$ combinations of the four local decision strategies in each of these three months), as well as extensions for the same problem up to seven monthly decision periods with the same numerical parameters.   

\begin{table}[htp!]
	\centering
	\caption{Results for the pig farm problem for different numbers of decision periods.}\label{tab:pigs_solutions}
	\setlength{\tabcolsep}{5.0pt}
	{\footnotesize{
\begin{tabular}{ccc}

\toprule

\# Months & Optimal value (DKK) & Solution time (s) \\ 
%% Note on new solutions
\midrule
3	& 764 & \p{00}0.01 \\
4	& 727 & \p{00}0.04 \\
5	& 703 & \p{00}0.62 \\
6	& 686 & \p{0}19.52  \\
7   & 674 & 617.21\\
%% Previous solutions
% 3	& 764 & \p{00}0.03 \\
% 4	& 727 & \p{00}0.09 \\
% 5	& 703 & \p{00}1.43 \\
% 6	& 686 & \p{0}43.83  \\
% 7   & 674 & 920.48\\
\bottomrule
\end{tabular}
}}
\end{table}

\begin{figure}
\centering
\includegraphics[width = 1.0\textwidth]{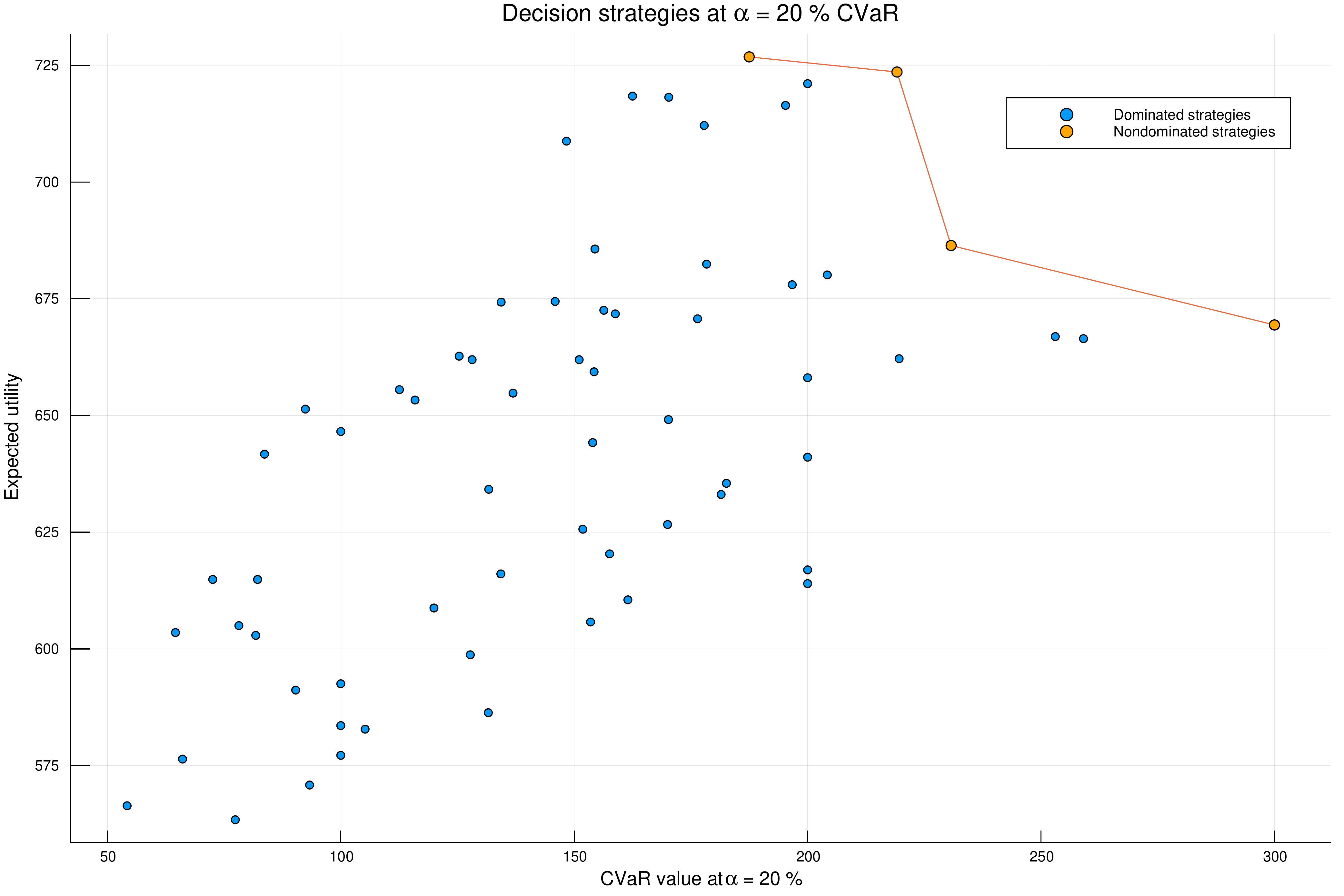}
\caption{Expected utilities and conditional expectations in the lower $\alpha = 0.20$ tail. The four non-dominated strategies are connected and marked with orange circles.}
\label{fig:cvar_pigs}
\end{figure}

Using the formulations in Section \ref{sec:multiple}, one can also determine the non-dominated strategies based on the consideration of the two objectives of maximizing (i) the overall expected utility and (ii) the conditional expectation in the lower $\alpha = 0.20$ tail. The 64 different strategies are presented in Figure \ref{fig:cvar_pigs} which shows the  corresponding overall expected utility (assuming risk-neutral preferences over monetary consequences) and the conditional CVaR expectation in the lower $\alpha = 20 \%$ tail of for each of the 64 decision strategies for this 4-month pig problem. However, the number of strategies grows rapidly with the number of periods, meaning that solving the problem through explicit enumeration becomes increasingly impractical. In the case of seven periods, for example, there are $4^7 = 16~384$ decision strategies.  

In Figure \ref{fig:cvar_pigs}, the four non-dominated strategies are connected and marked with orange circles, while the remaining 60 dominated strategies are marked with blue circles. Going from left to right, the first non-dominated strategy has the highest expected utility, while the fourth one has the highest conditional lower tail expectation. The vaccination policies in these non-dominated strategies are, respectively, as follows:
\begin{itemize}[itemsep=-1.0ex]
\item[1.] Never treat at 1st month. Treat at 2nd and 3rd month if and only if test results are positive.
\item[2.] Never treat at 1st and 2nd month. Treat at 3rd month if and only if test results are positive.
\item[3.] Never treat at 1st and 3rd month. Treat at 2nd month if and only if test results are positive.
\item[4.] Never treat at any of the 3 months.
%, regardless of test results.
\end{itemize}

Thus, the local strategy of never treating in the first month is a robust decision recommendation, because it is contained in all non-dominated strategies and consequently in the set of `core' selections in the Robust Portfolio Modelling (RPM) framework \citep{LiesioMildEtAl2007,LiesioMildEtAl2008}. Moreover, all local strategies involving treatments based on negative test results can be ruled out from consideration, because they are not in any non-dominated strategies and thus belong to the set of `exterior' RPM selections.    

\section{Summary and Conclusions}

In this paper, we have developed Decision Programming as an MILP optimization approach for solving mixed-integer multi-stage decision problems with discrete decisions and chance events. Such problems can be represented as influence diagrams, including LIMIDs in which the usual assumption of `no-forgetting' may not hold. In this approach, risk preferences can be captured through non-linear utility functions over consequences or, alternatively, by extending the objective function with terms for risk measures or by introducing risk constraints. Multiple objectives can be handled, for instance, by using a weighted additive linear function to aggregate consequences (or their utilities) across different value nodes. The set of all non-dominated strategies can also be computed with MILP by employing a weighted linear objective function together with the sequential introduction of constraints to eliminate dominated strategies as well as already discovered non-dominated strategies from further consideration.  

In the context of stochastic optimization, Decision Programming is particularly useful in mixed-integer decision problems where the probabilities in the scenario tree depend endogenously on earlier integer-valued decisions. This ability to handle endogenous uncertainties can be helpful, for instance, when appraising  R\&D and marketing investments, because the size of the market as well as the products' market performance are often contingent on these earlier decisions. From this perspective, the proposed approach can be viewed as a generalization of Contingent Portfolio Programming that allows the chance events describing the scenario tree depend on project selection decisions. 

Importantly, the  Decision Programming framework can be employed to address problems that cannot solved with dynamic programming techniques, such as problems in which earlier decisions cannot be recalled or in which the presence of deterministic and chance constraints make it impractical or impossible to apply dynamic programming techniques. Therefore, although Decision Programming has parallels to developments in stochastic mixed-integer dynamic programming (such as employing mathematical programming formulation to find optimal policies, as in the seminal work of \citet{Manne1960} and ensuing literature; see \citet{Bertsekas2012} for a thorough exposition), Decision Programming makes it possible to solve a broader class of problems which are not amenable to dynamic programming. Technically, the key feature of our approach is that, instead of exploiting recursion as the underpinning framework, we exploit the expressiveness of influence diagrams for problem structuring and then develop equivalent deterministic MILP formulations that can be solved using off-the-shelf MILP solvers.

Based on our numerical experiments, the Decision Programming approach allows for solving large-scale problems to optimality. Quite importantly, its computational performance can be radically enhanced through the use of probability cuts which exploit the specific properties of probabilistic constraints as well as whatever symmetric properties the problem structure may feature. 

Nevertheless, Decision Programming is subject to the well-known curse of dimensionality, just as other linear programming-based approaches for solving dynamic problems with a larger number of decision periods and uncertainties. However, given that powerful MILP decomposition and parallelization techniques are becoming widely accessible, the proposed approach holds considerable promise in extending the expressiveness of influence diagrams in problem structuring while offering possibilities for handling multiple objectives subject to a much broader range of constraints than what conventional approaches for building and solving influence diagrams can accommodate.

\section*{Acknowledgements} 

This research has been partly funded by the project {\it Platform Value Now} of the Strategic Council of the Academy of Finland (funding decision number 314207).

\appendix

% \section{Computational Environment}\label{app:CE}

% All the problem instances in Section \ref{sec:CT} were solved with an Intel Xeon E3-1230 v5 desktop clocked at 3.40 GHz with 32 GB RAM running Windows 10 x64 Education Edition. The Decision Programming approach was coded in Julia 1.1.0, and the mixed-integer linear programming problems were solved with Gurobi 8.1.0 using 2 out of 8 available threads.

 \section{Proofs}\label{app:proofs}

 \textsc{\noindent Proof of Theorem 1}. Let $Z \in \mathbb{Z}$ and take any path $s \in S$. The information set of the first node $k=1$ is empty. If this node is a chance node, the random variable $X_1$ does not depend on $Z$ and thus $\pi_1(s) = \mathbb{P}(X_1 = s_1) = \mathbb{P}(X_1 = s_1 \mid Z)$. If it is a decision node, there are two cases. First, if $Z_1 = Z_1(\emptyset) = s_1$, it follows that $\mathbb{P}(X_1=s_1 \mid Z) = 1$ while \eqref{zvariables} gives $z(s_1) = 1$. Thus, by \eqref{DecisionPi} we have $\pi_1(s) = z(s_1) = 1 = \mathbb{P}(X_1=s_1 \mid Z)$. Second, if  $Z_1 \neq s_1$, then $\mathbb{P}(X_1=s_1 \mid Z) = 0$ while $z(s_1) = 0$ gives $\pi_1(s) = 0$, and hence $\pi_1(s) = 0 = \mathbb{P}(X_1=s_1 \mid Z)$ in this case, too. Thus, Theorem \ref{scenariopaths} holds for $k=1$. 

 Assume that \eqref{pi_recursion} holds for $j \in 1,\dots,k-1$ with $k-1 < n$. We show that it holds for $k$, too. If $k \in C$ is a chance node, $\{ j \mid j \in D, j \leq k \} = \{ j \mid j \in D, j \leq k-1 \}$ and 
 \begin{eqnarray*}
 \mathbb{P}(s_{1:k} \mid Z) & = & \bigg( \prod_{\stackrel{i \in C}{i \leq k}} \mathbb{P} \big(X_i=s_i \mid X_{I(i)} = s_{I(i)} \big) \bigg) \bigg( \prod_{\stackrel{j \in D}{j \leq k}} \mathbb{I} \big( Z_j(s_{I(j)})=s_j \big) \bigg) \\
 & = & \mathbb{P} \big(X_k=s_k \mid X_{I(k)}=s_{I(k)} \big) \bigg( \prod_{\stackrel{i \in C}{i \leq k-1}} \mathbb{P} \big(X_i=s_i \mid X_{I(i)}=s_{I(i)} \big) \bigg) \bigg( \prod_{\stackrel{j \in D}{j \leq k-1}} \mathbb{I} \big( Z_j(s_{I(j)})=s_j \big) \bigg) \\
 & = & \mathbb{P} \big(X_k=s_k \mid X_{I(k)} = s_{I(k)} \big) \pi_{k-1}(s) = \pi_k(s), 
 \end{eqnarray*}
 where the last equality follows the induction hypothesis and  \eqref{condprob}. Analogously, if $k \in D$ is a decision node, then  
 \begin{eqnarray*}
 \mathbb{P}(s_{1:k} \mid Z) & = & \mathbb{I} \big( Z_k(s_{I(k)})=s_k \big) \bigg( \prod_{\stackrel{i \in C}{i \leq k-1}} \mathbb{P} \big(X_i=s_i \mid X_{I(i)}=s_{I(i)} \big) \bigg) \bigg( \prod_{\stackrel{j \in D}{j \leq k-1}} \mathbb{I} \big( Z_j(s_{I(j)})=s_j \big) \bigg) \\
 & = & z(s_k \mid s_{I(k)}) \pi_{k-1}(s) = \pi_k(s), 
 \end{eqnarray*}
 where the last equality follows the induction hypothesis and equations \eqref{zvariables} and \eqref{DecisionPi}. \halmos

\vskip1.0cm

{\noindent \sc Proof of Proposition 1}. Choose $\alpha \in (0,1]$ and consider $\eta^* = {\rm VaR}_{\alpha}(Z)$ which, by \eqref{VaR}, is well defined. Then constraints \eqref{beginCVaR} -- \eqref{CVaR4} are satisfied by $\rho(s), \overline{\rho}(s), \lambda(s)$ and $\overline{\lambda}(s)$, defined so that $\lambda(s) = \overline{\lambda}(s) = 1$ for paths such that $\mathcal{C}(s) < \eta^*$; $\overline{\lambda}(s) = 1$ and $\lambda(s) = 0$ for $\mathcal{C}(s) = \eta^*$; and $\lambda(s) = \overline{\lambda}(s) = 0$ for $\mathcal{C}(s) > \eta^*$. From \eqref{CVaR5}-\eqref{CVaR7} it follows that $\rho(s) = \overline{\rho}(s) = \pi(s)$ when $\mathcal{C}(s) < \eta^*$; $0 = \rho(s) = 0 \leq \overline{\rho}(s) \leq \pi(s)$ for $\mathcal{C}(s) = \eta^*$; and $\rho(s) = \overline{\rho}(s) = 0$ when $\mathcal{C}(s) > \eta^*$. By the choice of $\eta^*$, is it possible to choose variables $\overline{\rho}(s) \geq 0$ for $\mathcal{C}(s) = \eta^*$ so that \eqref{CVaR8} gives the correct tail expectation $\sum_{s \in S} \overline{\rho}(s)\mathcal{C}(s) / \alpha$ in \eqref{CVaR}. Finally, assume that there exists another solution for some $\eta' < \eta^*$. But then \eqref{CVaR8} implies that the probability $\alpha$ is attained as the sum of those paths whose consequence is lower than or equal to $\eta'$, violating the assumption that $\eta^* = {\rm VaR}_{\alpha}(Z)$. 
\halmos

\bibliography{references}

%%%%%%%%%%%%%%%%%
\end{document}